\documentclass[1p]{elsarticle}
\makeatletter
\def\ps@pprintTitle{%
	\let\@oddhead\@empty
	\let\@evenhead\@empty
	\def\@oddfoot{\centerline{\thepage}}%
	\let\@evenfoot\@empty}
\makeatother
\usepackage{moreverb}

\usepackage{graphicx} 
\usepackage{epsfig} 
\usepackage{amsmath} 
\usepackage{amssymb}  
\usepackage[scientific-notation=true]{siunitx}
\usepackage{commath}
\usepackage{xcolor} 
\usepackage{import}
\usepackage{epstopdf}
\usepackage{pgfplots}
\usepackage{tikz}
\usepackage{capt-of}
\usepackage{eurosym}
\usepackage{hyperref}
\usepackage{subcaption}
\graphicspath{{../fig/}{../imgs/}}
\newcommand\BibTeX{{\rmfamily B\kern-.05em \textsc{i\kern-.025em b}\kern-.08em
T\kern-.1667em\lower.7ex\hbox{E}\kern-.125emX}}

\begin{document}

\begin{frontmatter}
\title{An optimal hierarchical control scheme for smart 	generation units: an application to combined steam and electricity generation}
\fntext[t1]{Published on Journal of Process Control. DOI:\url{https://doi.org/10.1016/j.jprocont.2020.08.006}}
\author[1,2]{Stefano Spinelli\corref{cor1}}
\ead{stefano.spinelli@stiima.cnr.it}

\author[2]{Marcello Farina}

\author[1]{Andrea Ballarino}

\address[1]{Istituto di Sistemi e Tecnologie
		Industriali Intelligenti per il Manifatturiero Avanzato, Consiglio Nazionale delle Ricerche, Italy, Milano}

\address[2]{ Dipartimento di Elettronica, Informazione
		e Bioingegneria, Politecnico di Milano, Italy, Milano}
\cortext[cor1]{Corresponding author }

\begin{abstract}	
%
Optimal management of thermal and energy grids with fluctuating demand and prices 
requires to orchestrate the generation units (GU) among all their operating modes. A hierarchical approach is proposed to control  coupled energy nonlinear systems. The high level hybrid optimization defines the unit commitment, with the optimal transition strategy, and best production profiles. The low level dynamic model predictive control (MPC), receiving the set-points from the upper layer, safely governs the systems considering process constraints. 
To enhance the overall efficiency of the system, a method to optimal start-up the GU is here presented: a linear parameter varying MPC computes the optimal trajectory in closed-loop by iteratively linearising the system along the previous optimal solution. 
The introduction of an intermediate equilibrium state as additional decision variable permits the reduction of the optimization horizon,
while a terminal cost term steers the system to the target set-point. Simulation results show the effectiveness of the proposed approach.
\end{abstract}


\end{frontmatter}

\section{Introduction}\label{chap:Intro}
\subsection{Motivation}
The energy and utilities industry is facing a deep transformation, driven by environmental policies, characterized by the transition from centralized generation systems to distributed ones. This shift, supported by advanced digital technologies, requires the introduction of new paradigms for energy management and control. This includes a radical change in the user role, the integration of different energy vectors, and the change of focus for the definition of the control problems from the device level to a system level.\smallskip\\
First, decentralization of the energy production and the integration of small scale generation are promoting a paradigm shift in the consumers, which gain an active role both on the modification of their consumption patterns, i.e. in demand-response programs, and by producing and dispatching locally their energy.
In this new scenario, the consumer is asked to abandon the perspective of "Energy as a Commodity", where energy is considered always available and cost-effective on demand toward the paradigm of "Energy-as-a-Service". The latter is envisioned as one of the main strategies to reach the Circular Economy of the energy industries, together with industrial symbiosis - i.e. utilization of excess energy and side streams - and local-level cooperation \cite{Deloitte2018}.\\
This shift is allowed by moving the focus of managing, optimization, and control problems from the device-level to the system level, e.g., as in smart grids.
The generation nodes are composed of multiple integrated units that require coordination and control to fulfill the demand of the various consumers, which in turn may vary (in part) their demand based on optimization criteria.\\
This also promotes a scenario in which the production of various forms of utilities is more and more integrated. A step-ahead in system flexibility is built on a cross-sectorial integration of different energy vectors and the development of tools and technologies that enable efficient utilization of multi-dimensional  energy  systems \cite{SET-Plan2018}. For this reason, significant research efforts must be devoted to alternative energy carriers, e.g., compressed air, heating and cooling networks, and to their integration.
This allows, e.g., to extend the idea of smart grids to Smart Thermal Energy Grid (Smart-TEG).\smallskip\\
In this work we focus on a local instance of the Smart-TEG, where an industrial consumer, as in large companies and in industrial parks, includes an internal network of generation units (GU), which are responsible to supply energy - in various forms - to the production units.\\
Among the different configurations of a general thermal-energy grid, in this work we select a plant with steam and electricity flows. Steam is chosen as an high energy carrier since it is the most efficient media to transport thermal energy; it is extensively used in different sectors, e.g. chemical, medical textile and food. Besides, it is employed for the production of electrical power in Combined Cycle Power Plants (CCPP) and in co-generative systems.\\%
The integrated multi-utility configuration of the GU is here composed by a Fire-Tube Boiler (FTB), for the generation of steam, and by an Internal Combustion Engine (ICE), which operates as a combined heat and power system (CHP), providing both electrical and thermal energy. The selected case-study presents an interlink of the two generation systems: the thermal energy produced by the CHP can be diverted inside the boiler to support steam generation, through a controllable valve.\\
While the thermal and electricity demand defined by the scheduling of the production units must be fulfilled at any time, the mix between internally generated energy and net exchange with the national grid must be optimized, as well as the commitment of the generation units and their operating point.\\
In this context, the management and control of the integrated generation units become a necessary element to operate them safely and efficiently - both from economic and environmental perspective - as well as to balance the varying demand.\\
\subsection{State of the art}
In the field of thermal power generation, many research efforts have been devoted in the past decades both to solve regulation problems at a single device level and to address optimization and high level control with reference to integrated generation systems and cogeneration plants.\smallskip\\
Several works cover dynamic control problems related to steam systems. Nowadays, the deregulated electricity market and the integration of renewable sources in the generation mix impose dramatic fluctuations of the power demand even for these plants. Steam generators are nonlinear systems and the large variation of the demand profile imposes them to cover a wide range of operating points, and this requires the design of advanced control schemes for effective regulation.\\
Boilers, for example, are multi-input/multi-output (MIMO) systems for which control of both water level and pressure is required, respecting process constraints.
State of the art control is classically done by decoupled PID loops on water level and pressure, governed respectively by feedwater flow rate and fuel flow rate \cite{Dukelow2013}. Regarding drum boilers, note that the water level control is a severe issue due to the so called shrink-and-swell effect. Several works propose advanced strategies, e.g., \cite{Moon2009}, \cite{Lu2005} and references herein. \cite{Moon2009} proposes a Dynamic Matrix Control (DMC) approach based on step-response models to describe system dynamics.
To deal with the nonlinearities of the system, in \cite{Prasad2002} the authors present a hierarchical control structure with a high-level nonlinear model predictive control (NMPC) based on the nonlinear plant model \cite{Pedret2000}.\\
Other works, e.g., \cite{FB2005}, \cite{RV2008} and \cite{Hassanein2004}, address the control of fire-tube boilers (FTB) with similar control approaches but with a focus on the pressure loop, since in FTBs the level control is a minor problem because of the much larger free surface of water.\smallskip\\
High level optimization and control schemes of energy systems have also been studied. As discussed in \cite{Baldea14}, a strong integration of decision-making layers is advantageous, but also challenging, and therefore different approaches are available. For example, regarding cogeneration plants, consisting of a Boiler-Turbine (BT) configuration, \cite{Klauco2017} and \cite{Wu2014} propose multi-layer architectures based on reference governor algorithms.\\
For optimal management of thermal systems systems, it is often necessary to schedule activation, deactivation and transitional start-up/shut-down procedures, i.e., to define optimal Unit Commitment (UC) problems \cite{Baldick1995}, which is applied to different problems.\\
For example, for co-generative plants, in \cite{Mitra2013} the authors propose a direct Mixed-Integer Linear Programming (MILP) formulation, while in \cite{FT2004} an hybrid system approach is presented for a Combined Cycle Power Plants island composed of a gas and a steam turbine with a heat recovery steam generator. The work \cite{Ashok2002} presents a generalized formulation to determine the optimal operating strategy for any industrial cogeneration plant, for the equipment and load profiles of a typical petrochemical industry.\\
The work \cite{Hawkes2009} addresses the high-level UC of generators and storage within a micro-grid considering installation of a CHP, linking electricity and heat production, where the UC problem is formulated as a linear programming; it only considers control of the aggregate energy flows at the high level. Instead, in \cite{Anand2019} a particle swarm optimization method is used for the solution to a multi-objective UC problem of a CHP unit.\smallskip\\
Also optimization of start-up procedures is a major issue for guaranteeing optimal and safe system behaviour. While, in the past, shut-downs and start-ups of generation units, e.g. CCPP, were mainly imposed by maintenance actions and activation procedures were seldom performed, at present these plants are mostly under-utilized and go through several start-and-stop procedures to follow the demand. Optimizing the start-up procedure can lead to reducing the start-up time, minimizing the operative cost, avoiding thermal stress of the metal components, and reducing the environmental footprint, e.g., by limiting the fuel consumption and the emissions.\\
Since the calculation of the optimal trajectory for the entire procedure is computationally expensive for online implementation mainly due to the system non-linearity, many recent works propose offline optimization procedures. 
Model-based optimization is proposed for start-up procedures of coal-fired power plants~\cite{Hubel2017}, Heat Recovery Steam Generators~\cite{Sindareh2014} (HRSG), and turbines~\cite{Albanesi2006},~\cite{Casella2011}.
Nonlinear Programming (NLP) problems are proposed in ~\cite{Kru2004} and~\cite{Franke2003},and~\cite{Belkhir2015} for the open-loop optimization of a CCPP drum boiler, considering the thermal stress model. \cite{Taler2015} and~\cite{Taler2015_2} present an extremely detailed model of the thermal stresses for critical components, used for the off-line optimization of firing curves.\\
However, the approaches based on offline optimization cannot manage possible disturbances or compensate for modelling errors, possibly leading to safety and reliability concerns. On the other hand, closed loop approaches can assure disturbance rejection and compensate drifting or model mismatch. In \cite{Frank2009} a Nonlinear Model Predictive Control scheme is proposed, where the authors choose a prediction horizon that includes the whole start-up. In view of this, although this approach is very promising, computational complexity is still an issue.
\subsection{Paper contribution and structure}
In this paper a hierarchical scheme based on time scale separation is proposed to address the optimal management and the control of the generation units of a smart thermal energy grid, consisting of a FTB and a CHP. The hierarchical approach, similar to the one applied in \cite{Cominesi2017} to micro-grids, permits to tackle this structured problem, consisting of both UC and dynamic control.\\
 The main contribution of this research work is the extension of hierarchical approaches - mostly studied in the context of electrical systems - to steam generating units. To the best of our knowledge, this has not been previously investigated.
	 In addition, to comply with this aim, original modeling aspects related to boiler systems are presented, as well as specific assumptions introduced for layer integration. Moreover, the proposed method for optimizing the boiler startup, considering output measurements, is another novel contribution of this work.\\
The objective of UC is to schedule the discrete transitions between the different operating modes of the GU, based on a hybrid model of the system, to fulfill the utility demand and minimize the operative cost.\smallskip\\
Considering the result of the UC, the dynamic control algorithm aims at tracking the operating point trajectories generated by the upper layer and at regulating the system dynamics around these points rejecting disturbances. A model predictive control approach is proposed for this layer.
For low-level control we focus on the FTB system, which presents more challenging aspects than CHP regulation. CHP control, in fact, is just aimed principally at maintaining the rotational speed of the motor in presence of varying load.\smallskip\\
We also propose a nonlinear MPC approach for the control of the boiler start-up, inspired by the tracking MPC presented in \cite{Limon2008}. As it will be discussed, this allows to limit the optimization horizon, therefore reducing the numerical complexity of the related optimization problem.
To achieve an efficient implementation of the NMPC, we use here a parameter-varying linearisation of the nonlinear system, computed along the state/input trajectory obtained at the precedent optimization cycle.\smallskip\\
Preliminary results are proposed in \cite{Spinelli2018} and \cite{Spinelli2019}. In \cite{Spinelli2018} the hierarchical control structure for optimization of the GU behaviour is introduced. Differently from \cite{Spinelli2018}, in this paper the discrete high-level model includes an additional operating mode: the stand-by. An extensive analysis of process data has in fact highlighted the presence of a characteristic operating mode, distinguishable from both productive modes.\\
In \cite{Spinelli2019} the FTB start-up optimization method discussed here is introduced, for a simpler boiler model and without integration in the comprehensive hierarchical control scheme. In this paper we take several steps forward with respect to these preliminary works, also as far as the models used for FTB dynamic control are concerned. In fact, here we
consider an enhanced dynamical model of the system, in which the parameters are identified on available experimental data. Also, the FTB dynamical model is extended with the inclusion of the steam header model, required to represent the actual configuration of the system and the sensor location, and hence improving the description and reconstruction of available measurements. Furthermore, the problem has been addressed in a more realistic set-up: the state feedback assumption has been removed and an observer has been set up considering the available sensors.\\
The paper is organized as follows: Section \ref{chap:model} describes in detail the model of the GU, presenting the subsystem models used at both low and high hierarchical levels,  while Section \ref{chap:hier_ctrl} describes in detail the proposed control structure. Section \ref{chap:ltv_mpc} presents the optimal predictive control method used in the lower level for start-up optimization, while in Section \ref{chap:sim_res} simulation results are reported.
Finally, in Section \ref{chap:concl} some conclusions are drawn.

\section{Models of the generation units}
\label{chap:model}
The GU case study considered in this work is composed of a CHP and a FTB physically connected together. Specifically, the hot exhaust gases of the ICE can be either expelled in the air or directed inside the boiler thanks to a commanded valve. The residual thermal energy of the flue gases can thus be recovered in the FTB to sustain the steam production with a lower primary input consumption.\\
The mathematical models of these systems are discussed in the following sections: first, state-space models, derived from first principles and considered in the lower control layer, are presented for the CHP and the FTB. Secondly, the hybrid models \cite{Lygeros2008}, used  at UC level, are described.\smallskip\\ %
Different sampling times are used at the two different layers: the dynamic models used at low level are discretized with a sample time of $\tau_{\rm\scriptscriptstyle L}=60$ s, while at high-level the sampling time is $\tau_{\rm\scriptscriptstyle H}=15$ min. 
When necessary for better clarity, we will distinguish continuous-time from discrete-time variables using the following notation. For a generic variable $v$, $v(t)$ denotes the value of $v$ at time instant $t\in\mathbb{R}$, while $v[h]$ (respectively $v[k]$) denotes the value of $v$ at low-level discrete time step $h\in \mathbb{Z}$ (respectively high-level discrete time step $k\in \mathbb{Z}$).
\subsection{The CHP dynamic model in full operation mode}
\label{sec:Model-sub:CHP_dyn}
\begin{figure}[thbp]
	\centering
	\includegraphics[width=1\linewidth]{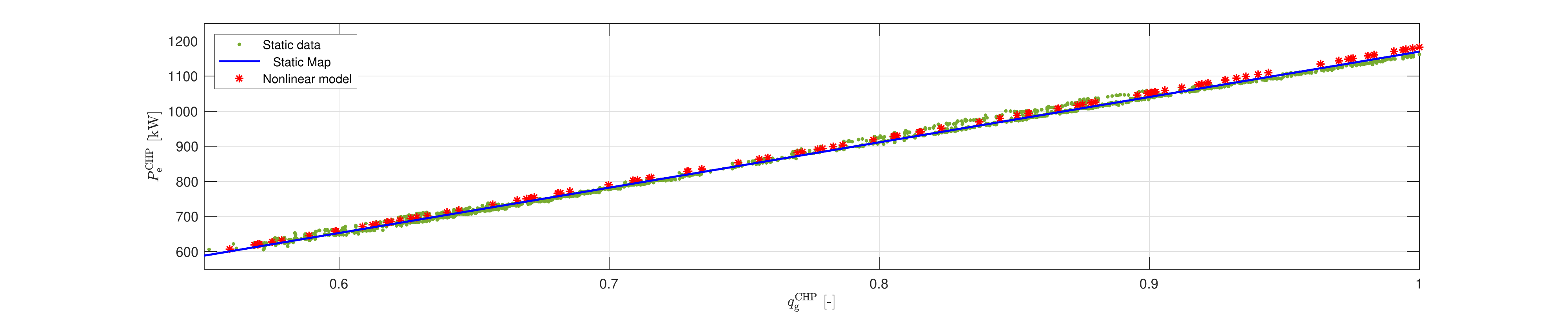}
	\caption{Relationship between $q_{\rm \scriptscriptstyle g}^{\rm \scriptscriptstyle CHP}$ and $P_{\rm \scriptscriptstyle e}^{\rm \scriptscriptstyle CHP}$ in steady-state conditions. Green dots: experimental data; red dots: static map of the nonlinear model devised in Section~\ref{sec:Model-sub:CHP_dyn}; blue line: identified affine relationship \eqref{eq:IN-OUT-CHP}. $q_{\rm \scriptscriptstyle g}^{\rm \scriptscriptstyle CHP}$ is adimensionalized for confidentiality.}
	\label{fig:Static_CHP}
\end{figure}

The CHP is a natural gas Internal Combustion Engine (ICE) connected to a three phase synchronous generator, to provide electricity.
The considered model of the ICE in operation mode is derived from the works of \cite{Guzzella2004} and \cite{Videla2007}: a lumped parameter model of a internal combustion engine, i.e. a mean value model, is connected to a static alternator.
The dynamical model of the ICE considers two states, the angular speed $\omega$ of the motor and the engine internal pressure $p^{\rm \scriptscriptstyle CHP}$.
The inputs are the position of the throttle valve, which is assumed to be directly related to the natural gas mass flow rate $q^{\rm \scriptscriptstyle CHP}_{\rm \scriptscriptstyle g}$ by stoichiometric combustion assumption, and the electrical demand $P_{\rm \scriptscriptstyle e}^{\rm \scriptscriptstyle CHP}$.
In addition to the state variables, another output channel is considered: the exhaust gas enthalpy $H_{\rm\scriptscriptstyle ex}^{\rm\scriptscriptstyle CHP}$ that depends on the states of the system.\\
Stoichiometric ratio, volumetric and thermal efficiency are approximated by second order polynomials of the systems states. Similarly, the heat transferred from the combustion gases to the ICE walls, used to evaluate the remaining enthalpy of the exhaust gases is simplified as a polynomial map considering motor angular speed and its internal pressure.\\
System parameters are derived from the physical and geometrical properties of the CHP, starting from literature values, as in~\cite{Guzzella2004}. Moreover, the gain matrix has been fine tuned based on available static data: in Figure \ref{fig:Static_CHP} the static map obtained by available data is compared to the corresponding map computed using the identified model. Values on the $x$ axis are adimensional for confidentiality reasons.\\
The local CHP controller, in this paper, is a standard PI controller acting on the throttle valve input, i.e. equivalently to the natural gas flow rate $q_{\rm \scriptscriptstyle g}^{\rm \scriptscriptstyle CHP}$, to regulate the angular velocity $\omega$ to its nominal value $\bar{\omega}= 50 Hz$. The PI loop is tuned to obtain a settling time is smaller than $1$~s.\\
Thanks to this regulating loop, the CHP is characterized by one input only, i.e., the electrical demand $P_{\rm \scriptscriptstyle e}^{\rm \scriptscriptstyle CHP}$, while the gas flow rate $q_{\rm \scriptscriptstyle g}^{\rm \scriptscriptstyle CHP}$ and the exhaust gas enthalpy $H_{\rm\scriptscriptstyle ex}^{\rm\scriptscriptstyle CHP}$ are to be considered as system outputs.\\
It is worth noting that the settling time of the system, $0.5$~s about, is much faster than the sampling time, $\tau_{\rm\scriptscriptstyle L}$, of the ''fast'' MPC controller designed in the following (see Section~\ref{sec:control-sub:LL}).
Based on the MPC timescale, we can assume the ICE to be always in steady state conditions, therefore we can model the CHP behaviour just considering the static gains where the power produced by the ICE is equal to $P_{\rm \scriptscriptstyle e}^{\rm \scriptscriptstyle CHP}$.\\
In addition, the system in full operation mode is constrained to work within a range of produced power; as a consequence, the input $q_{\rm \scriptscriptstyle g}^{\rm \scriptscriptstyle CHP}$ is constrained to lie in a set, where also a lower bound is defined.
\begin{equation}
q_{\rm \scriptscriptstyle g}^{\rm \scriptscriptstyle CHP}\in[\bar{q}_{\rm \scriptscriptstyle g \, min}^{\rm \scriptscriptstyle CHP},\bar{q}_{\rm \scriptscriptstyle g \, max}^{\rm \scriptscriptstyle CHP}]
\label{eq:CHP_CONSTRAINT}
\end{equation}
where $\bar{q}_{\rm \scriptscriptstyle g \, max}^{\rm \scriptscriptstyle CHP}>\bar{q}_{\rm \scriptscriptstyle g \, min}^{\rm \scriptscriptstyle CHP}>0$.

\subsection{The boiler dynamic model} \label{sub:boiler_low}
The mathematical model of the fire tube boiler, used for the lower level control of the system and inspired by~\cite{Ortiz2011} and~\cite{Astrom2000}, is described in the following sections.
\subsubsection{The boiler physical equations}
 In the peculiar configuration of the case study, because of the interconnection between the CHP and the boiler, we consider here two sources of heat, i.e., the gas burner and the exhaust gases diverted from the CHP. The set of equations - based on mass and energy conservation laws - describing the system dynamics is
\begin{subequations}
\begin{align}
\frac{d}{dt}  [\rho_{\scriptscriptstyle \rm s} V_{\scriptscriptstyle \rm s} + \rho_{\scriptscriptstyle \rm w} V_{\scriptscriptstyle \rm w}]  = & q_{\scriptscriptstyle \rm f} - q_{\scriptscriptstyle \rm s}^{\scriptscriptstyle \rm B} \label{eqn:cons_eq_in1}\\
\frac{d}{dt}  [\rho_{\scriptscriptstyle \rm s} e_{\scriptscriptstyle \rm s} V_{\scriptscriptstyle \rm s} + \rho_{\scriptscriptstyle \rm w} e_{\scriptscriptstyle \rm w} V_{\scriptscriptstyle \rm w} ]  = & Q^{\rm\scriptscriptstyle t,GAS}_{\scriptscriptstyle \rm{m} \rightarrow \rm{w}} + Q^{\rm\scriptscriptstyle t,CHP}_{\scriptscriptstyle \rm{m} \rightarrow \rm{w}} + q_{\scriptscriptstyle \rm f} h_{\scriptscriptstyle \rm f} - q_{\scriptscriptstyle \rm s}^{\scriptscriptstyle \rm B} h_{\scriptscriptstyle \rm s} \label{eqn:cons_eq_in2}\\
\frac{d}{dt}  [M^{\rm\scriptscriptstyle t,GAS} C_{\rm\scriptscriptstyle p} T^{\rm\scriptscriptstyle t,GAS}]  = & Q^{\rm\scriptscriptstyle B}_{\rm\scriptscriptstyle g} - Q^{\rm\scriptscriptstyle t,GAS}_{\scriptscriptstyle \rm{m} \rightarrow \rm{w}} \label{eqn:cons_eq_in3}\\
\frac{d}{dt}  [M^{\rm\scriptscriptstyle t,CHP} C_{\rm\scriptscriptstyle p} T^{\rm\scriptscriptstyle t,CHP}]  = & u_{\rm\scriptscriptstyle ex} \eta_{\scriptscriptstyle \rm ex} H^{\scriptscriptstyle \rm CHP}_{\scriptscriptstyle \rm ex} - Q^{\rm\scriptscriptstyle t,CHP}_{\scriptscriptstyle \rm{m} \rightarrow \rm{w}}\label{eqn:cons_eq_in4}
\end{align}
\label{eqn:cons_eq_in}
\end{subequations}
where $\rho$ is the density, $V$ the volume, $q$ the mass flow rate, $T$ the temperature, $h$ the specific enthalpy and $e$ the internal energy, of each system component denoted by its relative subscript notation: gas ($\rm{g}$), tubes ($\rm{t}$), steam ($\rm{s}$), water ($\rm{w}$), and feed-water ($\rm{f}$).\\ Equations \eqref{eqn:cons_eq_in3}-\eqref{eqn:cons_eq_in4} represent the energy balance of the heating tubes, where $m_{\rm\scriptscriptstyle t,GAS}$ (respectively, $m_{\rm\scriptscriptstyle t,CHP}$) denotes the mass of the tubes heated by the gas burner (respectively, transporting exhaust gases from the CHP), while $T^{\rm\scriptscriptstyle t,GAS}$ (respectively, $T^{\rm\scriptscriptstyle t,CHP}$) denotes their temperature. Moreover, the status of the interlink valve (open/closed), which can divert the flue gases from the CHP to the FTB is given by $u_{\rm\scriptscriptstyle ex} \in \{0, 1\}$ while $\eta_{\scriptscriptstyle \rm ex}$ denotes the rate of enthalpy of the exhaust gases that can be transmitted to the tubes.\\ The heat power of the gas burner is related to the natural gas flow rate $q^{\scriptscriptstyle \rm B}_{\rm\scriptscriptstyle g}$ by the following relation
$$Q^{\rm\scriptscriptstyle B}_{\rm\scriptscriptstyle g} = \eta^{\rm\scriptscriptstyle B} c_{\rm\scriptscriptstyle LHV} q^{\scriptscriptstyle \rm B}_{\rm\scriptscriptstyle g}$$
where $c_{\rm\scriptscriptstyle LHV}$ is the lower heating value. The heat transmitted from the metal walls to the water can be modelled as
\begin{align*}
Q^{\rm\scriptscriptstyle t,GAS}_{\scriptscriptstyle \rm{m} \rightarrow \rm{w}} &= \beta (T^{\rm\scriptscriptstyle t,GAS} - T_{\scriptscriptstyle \rm{w}})\\
Q^{\rm\scriptscriptstyle t,CHP}_{\scriptscriptstyle \rm{m} \rightarrow \rm{w}} &= \beta (T^{\rm\scriptscriptstyle t,CHP}-T_{\scriptscriptstyle \rm{w}})
\end{align*}  where $\beta$ is the heat transfer coefficient depending on the boiling two phase mixture of steam and water, close to the tube walls, that induces a natural recirculation and the interaction of numerous tubes in the bundle, based on Cooper correlation~\cite{Bejan2003}.\\
The conservation laws \eqref{eqn:cons_eq_in} can be simplified and further manipulated considering in the FTB two separated regions for water and steam and assuming the thermodynamic equilibrium between the regions.
It is thus possible to split the equation \eqref{eqn:cons_eq_in1} for the two regions:
\begin{align}
\frac{d}{dt}  [\rho_{\scriptscriptstyle \rm s} V_{\scriptscriptstyle \rm s} ]  = & q_{\scriptscriptstyle \rm s}^{ \rm\scriptscriptstyle  w \rightarrow s} - q_{\scriptscriptstyle \rm s }^{\scriptscriptstyle \rm B}\\
\frac{d}{dt}  [ \rho_{\scriptscriptstyle \rm w} V_{\scriptscriptstyle \rm w}]  = & q_{\scriptscriptstyle \rm f} - q_{\scriptscriptstyle \rm s} ^{ \rm\scriptscriptstyle   w \rightarrow s} \label{eq:mass_w}
\end{align}
where $q_{\scriptscriptstyle \rm s}^{\rm\scriptscriptstyle w \rightarrow s}$ is the steam released from the water region into the steam section.\\ By neglecting steam accumulation in the steam zone $q_{\scriptscriptstyle \rm s}^{\rm\scriptscriptstyle w \rightarrow s} = q_{\scriptscriptstyle \rm s } $ and by writing the water mass balance as:
\begin{equation*}
\frac{d}{dt}  [ \rho_{\scriptscriptstyle \rm w} V_{\scriptscriptstyle \rm w}]  = \rho_{\scriptscriptstyle \rm w}\frac{d}{dt}  V_{\scriptscriptstyle \rm w} + V_{\scriptscriptstyle \rm w} \frac{d}{dt}  \rho_{\scriptscriptstyle \rm w}  \simeq - \rho_{\scriptscriptstyle \rm w}\frac{d V_{\scriptscriptstyle \rm s}}{dl_{\scriptscriptstyle \rm w}}  \frac{dl_{\scriptscriptstyle \rm w}}{dt}
\end{equation*}
we can simplify the formulation.
The above approximation can be done assuming the second term negligible and remembering that $V_{\scriptscriptstyle \rm w} + V_{\scriptscriptstyle \rm s} = V_{\scriptscriptstyle \rm Tot}$ is constant
\begin{equation*}
\frac{d V_{\scriptscriptstyle \rm w}}{dt} = - \frac{d V_{\scriptscriptstyle \rm s}}{dt}
\end{equation*}
Based on FTB geometry, is it possible to express the volume of the steam zone $V_s$ as a function of the water level, $l_{\scriptscriptstyle \rm w}$:
\begin{equation*}
V_s = A_s L \simeq \frac{4+\pi}{2}r^2L - 2rLl_{\scriptscriptstyle \rm w}
\end{equation*}
The same separation in steam and water regions can be considered also for the energy balance in eq.~\eqref{eqn:cons_eq_in2}. Recalling the assumption of thermal equilibrium  for the two regions and of saturated liquid, $T_{\scriptscriptstyle  w} = T_{\scriptscriptstyle  s}$, we can write the energy balance related only to the water region:
\begin{equation*}
\frac{d}{dt}  [\rho_{\scriptscriptstyle \rm w} V_{\scriptscriptstyle \rm w} e_{\scriptscriptstyle \rm w}] = M_{\scriptscriptstyle \rm w} \frac{d}{dt}  [ e_{\scriptscriptstyle \rm w}]
+e_{\scriptscriptstyle \rm w} \frac{d}{dt}  [\rho_{\scriptscriptstyle \rm w} V_{\scriptscriptstyle \rm w} ]
\end{equation*}
The above energy conservation law for the water side can be expressed in terms of water temperature by recalling eq.~\eqref{eq:mass_w} and the thermodynamic relations for enthalpy, $dh = c_{\rm\scriptscriptstyle p} dT$, and internal energy, $de = c_{\scriptscriptstyle \rm v}dT$.\\
The nonlinear dynamic model of the FTB boiler can be recast in the following form:
\begin{subequations}
\label{eqn:cons_eq_in_bis}
\begin{align}
\frac{dT^{\rm\scriptscriptstyle t,GAS}}{dt}   = &\frac{1}{M^{\rm\scriptscriptstyle t,GAS} c_{\rm\scriptscriptstyle p} } [Q^{\rm\scriptscriptstyle B}_{\rm\scriptscriptstyle g} - \beta (T^{\rm\scriptscriptstyle t,GAS} - T_{\scriptscriptstyle \rm{w}})]\\
\frac{dl_{\rm\scriptscriptstyle w}}{dt}  = &\frac{1}{2rL} [q_{\scriptscriptstyle \rm f} - q_{\scriptscriptstyle \rm s}^{\scriptscriptstyle \rm B}]\label{eqn:cons_eq_in_bis2}\\
\frac{dT_{\rm\scriptscriptstyle w}}{dt}    =  &\frac{1}{\rho_{\rm\scriptscriptstyle w} V_{\rm\scriptscriptstyle w}(l_{\rm\scriptscriptstyle w}) c_{\rm\scriptscriptstyle v}} [Q^{\rm\scriptscriptstyle t,GAS}_{\scriptscriptstyle \rm{m} \rightarrow \rm{w}} + Q^{\rm\scriptscriptstyle t,CHP}_{\scriptscriptstyle \rm{m} \rightarrow \rm{w}}
+ q_{\scriptscriptstyle \rm f} c_{\rm\scriptscriptstyle p} (T_{\rm\scriptscriptstyle f} - T_{\rm\scriptscriptstyle w}) - q_{\scriptscriptstyle \rm s}^{\scriptscriptstyle \rm B} h_{\rm\scriptscriptstyle s}(T_{\rm\scriptscriptstyle w})] \label{eqn:cons_eq_in_bis3}\\
\frac{dT^{\rm\scriptscriptstyle t,CHP}}{dt}  = &\frac{1}{M^{\rm\scriptscriptstyle t,CHP} c_{\rm\scriptscriptstyle p} } [u_{\scriptscriptstyle \rm ex} \eta_{\scriptscriptstyle \rm ex} H^{\scriptscriptstyle \rm CHP}_{\scriptscriptstyle \rm ex} - \beta (T^{\rm\scriptscriptstyle t,CHP} - T_{\scriptscriptstyle \rm{w}})]
\end{align}
\end{subequations}
The FTB is equipped with a set of sensors, collecting the steam pressure $p_{\rm\scriptscriptstyle s}$ and the water level $l_{\rm\scriptscriptstyle w}$, assumed to be available by the internal controller. The output equation for the pressure can be defined by considering that the system works at saturated conditions, therefore, the water and steam pressures satisfy the equality:
\begin{equation}
p_{\rm\scriptscriptstyle s} = p_{\rm\scriptscriptstyle w}=p_{\rm\scriptscriptstyle sat }(T_{\rm\scriptscriptstyle w}) \label{eq:out_eq1}
\end{equation}
where the function of the thermodynamic properties of the steam and water at saturation, $p_{\rm\scriptscriptstyle sat }(T)$, is defined  based on the Industrial Formulation IAPWS-IF97~\cite{SteamTab}.

The state vector of the complete FTB model, 
$x=[T^{\rm\scriptscriptstyle t,GAS}, l_{\rm\scriptscriptstyle w}, T_{\rm\scriptscriptstyle w}, T^{\rm\scriptscriptstyle t,CHP}]'$, includes the temperature of the heating tubes, the water level and its temperature.
The flow rates of feed-water, $q_{\scriptscriptstyle \rm f}$, natural gas,  $q^{\scriptscriptstyle \rm B}_{\rm\scriptscriptstyle g}$, and steam output, $q_{\rm\scriptscriptstyle s }^{\scriptscriptstyle \rm B}$, are the manipulable inputs of the model. It has to be noticed that, based on the operating mode of the system, the latter input can be considered either an actual  manipulated variable, e.g. as in start-up modes, or a measured disturbance as in full operating mode, when steam outflow is defined by the consumer demand.\\
The controlled output vector is composed by internal steam pressure and water level, i.e., $y = [l_{\scriptscriptstyle \rm w}, p^{\scriptscriptstyle \rm B}_{\scriptscriptstyle \rm s }]$, where the boiler steam pressure is computed by the nonlinear thermodynamic map at saturation \eqref{eq:out_eq1}.
The continuous time nonlinear model of the fire-tube boiler is, in short:
\begin{subequations}
\label{eq:nl_model_tot}
	\begin{align}
	\dot{x} &= f(x, u, \theta) \label{eq:nl_model}\\
	y &= 	g (x, u, \theta) \label{eq:nl_model_out}
	\end{align}
\end{subequations}%
where $\theta$ is a vector collecting the uncertain system parameters, specified in the following in Section~\ref{subsub:par_id}. Input and state variables are all subject to constraints:
\begin{subequations}
	\label{eq:boiler_CONSTRAINT}
	\begin{align}
	q_{\scriptscriptstyle \rm f}&\in [0,\bar{q}_{\scriptscriptstyle \rm f\,max}]\\
	q^{\scriptscriptstyle \rm B}_{\rm\scriptscriptstyle g}&\in [\bar{q}^{\scriptscriptstyle \rm B}_{\scriptscriptstyle \rm g \, min}, \bar{q}^{\scriptscriptstyle \rm B}_{\scriptscriptstyle \rm g \, max}]\\
	q_{\scriptscriptstyle \rm s}^{\scriptscriptstyle \rm B}&\in [\bar{q}_{\scriptscriptstyle \rm s \,min},\bar{q}_{\scriptscriptstyle \rm s \,max}] \cup \{0\}\\
	p_{\scriptscriptstyle \rm s}^{\scriptscriptstyle \rm B}&\in [\bar{p}_{\scriptscriptstyle \rm min}, \bar{p}_{\scriptscriptstyle \rm max}] \label{eq:p_const}\\
	l_{\scriptscriptstyle \rm w}&\in [\bar{l}_{\scriptscriptstyle \rm min}, \bar{l}_{\scriptscriptstyle \rm max}] \label{eq:l_const}
	\end{align}
\end{subequations}%
For control design purposes, the nonlinear dynamical model~\eqref{eq:nl_model} is discretized, with sampling time $\tau_{\rm L}$, using the forth-order Runge-Kutta method. With some abuse of notation, the so-obtained discrete-time model is
\begin{equation}
\label{eq:nl_RK4}
x[h+1] = f_{RK4}(x[h],u[h],\theta)
\end{equation}
where $h$ represents the discrete time step.
\subsubsection{Parameter identification} \label{subsub:par_id}
\begin{figure}[tbh]
	\centering
	\includegraphics[width=.7\linewidth]{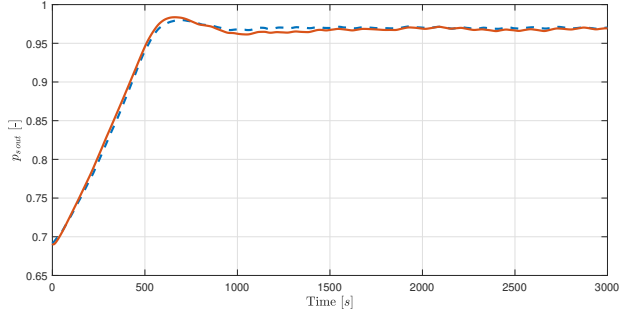}
	\caption{Comparison of the simulated steam header pressure $\tilde{y}$ (solid red) and experimental data $\tilde{y}_{\rm\scriptscriptstyle meas}$ (dashed blue)}
	\label{fig:Par_Id}
\end{figure}

The parameters of the model are computed based on the physical and geometric properties of the system. However, a fine tuning of some key parameters is conducted based on available data.\\
The data available at present for parameter identification purposes, however, do not correspond with all inputs and outputs of model \eqref{eq:nl_model_tot}: in fact, the pressure $p_{\scriptscriptstyle \rm s}$ and the water level $l_{\scriptscriptstyle \rm w}$ measurements, as well as the steam flow rate $q_{\scriptscriptstyle \rm s}^{\scriptscriptstyle \rm B}$ are not logged and available for identification. Nevertheless, in the steam header H, physically placed downstream of the FTB, a steam pressure $p_{\scriptscriptstyle \rm s}^{\scriptscriptstyle \rm H}$ and a flow rate $q_{\scriptscriptstyle \rm s}^{\scriptscriptstyle \rm H}$ sensors are present and connected to the plant data acquisition and monitoring system.
As a preliminary step to parameter identification, therefore, we first derive the model, linking the state-space model~\eqref{eq:nl_model_tot} variables with $p_{\scriptscriptstyle \rm s}^{\scriptscriptstyle \rm H}$ and $q_{\scriptscriptstyle \rm s}^{\scriptscriptstyle \rm H}$.
\paragraph{Model adaptation for parameter identification}\hfill\\
The steam header is a conduct, modelled as a chamber acting as a small accumulator. Here the mass conservation law and the compressible Bernoulli equation hold, i.e.,
\begin{align}
\frac{dM}{dt} &= q_{\rm\scriptscriptstyle s }^{\scriptscriptstyle \rm B} - q_{\rm\scriptscriptstyle s }^{\scriptscriptstyle \rm H} \label{eq:Cond_mass}\\
\frac{1}{2}\left( \frac{q_{\rm\scriptscriptstyle s }^{\scriptscriptstyle \rm B}}{\rho_{\rm\scriptscriptstyle s }^{\scriptscriptstyle \rm B} A^{\scriptscriptstyle \rm B}}\right)^2 + h_{\rm\scriptscriptstyle s }^{\scriptscriptstyle \rm B} + \Psi_{\rm\scriptscriptstyle s }^{\scriptscriptstyle \rm B} &= \frac{1}{2}\left( \frac{q_{\rm\scriptscriptstyle s }^{\scriptscriptstyle \rm H}}{\rho_{\rm\scriptscriptstyle s }^{\scriptscriptstyle \rm H} A^{\scriptscriptstyle \rm H}}\right)^2 + h^{\scriptscriptstyle \rm H}_{\rm\scriptscriptstyle s} + \Psi^{\scriptscriptstyle \rm H}_{\rm\scriptscriptstyle s} \label{eq:Comp_Bern}
\end{align}
where the superscripts $^{\rm{B}}$ and $^{\rm{H}}$ refer to the boiler outflow and the sensor location in the steam header. $M$ is the steam mass accumulated in the conduct, $A$ is the conduct area, $h$ is the enthlapy and $\Psi$ the potential associated with the conservative forces, i.e. the gravitational fields, which contribution is negligible in this case, as constant.\\
First, eq.~\eqref{eq:Cond_mass}, under the isochoric assumption, allows to write the expression of $q_{\rm\scriptscriptstyle s }^{\scriptscriptstyle \rm B}$ in \eqref{eqn:cons_eq_in_bis2} and \eqref{eqn:cons_eq_in_bis3} as a function of $q_{\rm\scriptscriptstyle s }^{\scriptscriptstyle \rm H}$ and the boiler variables as follows
\begin{equation}
q_{\rm\scriptscriptstyle s }^{\scriptscriptstyle \rm B} = V_{\rm\scriptscriptstyle c} \frac{d\rho}{dt} + q_{\rm\scriptscriptstyle s }^{\scriptscriptstyle \rm H} \label{eqn:q_in}
\end{equation}
where $V_c$ is the conduct volume and the density time-derivative at the boiler exit is computed as
\begin{equation}
\left. \frac{d\rho}{dt} \right|_{\rm\scriptscriptstyle s \,B}  =  \left.  \frac{d\rho}{dT} \right|_{\rm\scriptscriptstyle s \,B}  \left. \frac{dT}{dt} \right|_{\rm\scriptscriptstyle s \,B}\label{eqn:q_in2}
\end{equation}
After replacing~\eqref{eqn:q_in} and~\eqref{eqn:q_in2} in~\eqref{eqn:cons_eq_in_bis}, the resulting dynamic system is indicated as follows for brevity
\begin{equation}
    \dot{x} = \tilde{f}(x, \tilde{u}, \theta) \label{eq:nl_model_alt}
\end{equation}
where $\tilde{u}= [q_{\rm\scriptscriptstyle f }, q_{\rm\scriptscriptstyle g }^{\scriptscriptstyle \rm B}, q_{\rm\scriptscriptstyle s }^{\scriptscriptstyle \rm H}]$ is the ''new'' input vector.
Secondly, by recalling that enthalpy is a function of the pressure $h = h(p)$, 
we obtain, from~\eqref{eq:Comp_Bern}, the following Bernoulli equation
\begin{equation}
\frac{1}{2}\left( \frac{q_{\rm\scriptscriptstyle s }^{\scriptscriptstyle \rm B}}{\rho_{\rm\scriptscriptstyle s }^{\scriptscriptstyle \rm B}(p_{\rm\scriptscriptstyle s }^{\scriptscriptstyle \rm B}) A^{\scriptscriptstyle \rm B}}\right)^2 + h^{\scriptscriptstyle \rm B}_{\rm\scriptscriptstyle s}(p_{\rm\scriptscriptstyle s }^{\scriptscriptstyle \rm B}) = \frac{1}{2}\left( \frac{q_{\rm\scriptscriptstyle s }^{\scriptscriptstyle \rm H}}{\rho_{\rm\scriptscriptstyle s }^{\scriptscriptstyle \rm H}(p_{\rm\scriptscriptstyle s }^{\scriptscriptstyle \rm H}) A^{\scriptscriptstyle \rm H}}\right)^2 + h^{\scriptscriptstyle \rm H}_{\rm\scriptscriptstyle s}(p_{\rm\scriptscriptstyle s }^{\scriptscriptstyle \rm H}) \label{eq:Comp_Bern_out}
\end{equation}
Equation~\eqref{eq:Comp_Bern_out}, in turn, allows to derive the explicit expression of the header steam pressure $p_{\rm\scriptscriptstyle s }^{\scriptscriptstyle \rm H}$, regarded as a new output, as a function of the boiler system~\eqref{eq:nl_model_alt} variables $x$ and $\tilde{u}$. In short, we write $p_{\rm\scriptscriptstyle s }^{\scriptscriptstyle \rm H}=\tilde{y} = \tilde{g} (x, \tilde{u}, \theta)$.\\
\paragraph{Identification procedure}\hfill\\
The data at our disposal for identification purposes cover several days of plant operation, logged by the data acquisition and monitoring system in a historical database. In addition to the steam pressure and flow rate sensors located in the steam header, a flow-rate and a temperature sensors are available on the  feed-water conduct, as well as a flow-rate sensor located downstream of the junction from the main natural gas pipe: these sensors are used for the identification of the parameters of the boiler model. Moreover, since the variation of the feed water temperature is very limited, less than $1\%$, the average value is used in the boiler model as a fixed parameter.
To this aim, we selected data-sets related to the last part of the boiler start-up. 
The system, in this particular process condition, is more excited in terms of input variation and also the signal to noise ratio is larger, since the output of the system varies significantly. \\
 Moreover, the startup is characterized by a first part operated in open-loop, as the activation procedure is initially manual, and by a second half in closed-loop, as the regulator is triggered when water temperature reaches the regime value.
	The manual procedure, depicted in Figure \ref{fig:man_in}, reveals that the gas flow-rate is set to  three specific levels based on the water temperature. 
	The system operates in open-loop until $T_{\rm w}=0.95T_{\rm w}^{\rm SP}$, that triggers the activation of the PI regulator.\\
	Therefore, the dataset used for identification includes both data collected in open loop and in closed loop operation.\\
The  data-set used for training consists of arrays of $N=900$ samples, sampled with a sampling time of $\tau_{\rm\scriptscriptstyle id}=10$ s, of variables $\tilde{y}$ and $\tilde{u}$. We denote such data as $\tilde{u}_{\rm\scriptscriptstyle meas}[h]$ and $\tilde{y}_{\rm\scriptscriptstyle meas}[h]$, where $h=1,2,\dots,N$.\\
\begin{figure}[bht]
	\centering
		\includegraphics[width=0.6\linewidth]{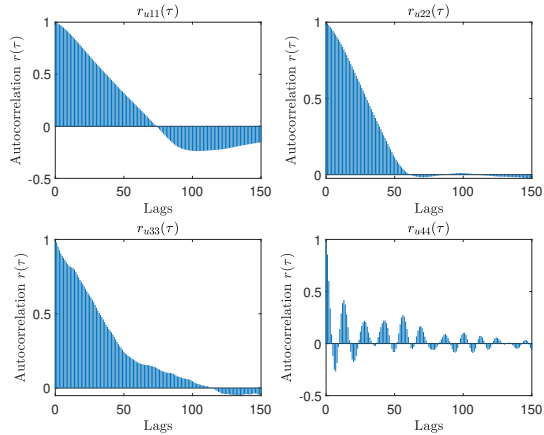}
		\caption{Spectral density of the input data used in training for the identification of model parameters}
		\label{fig:Par_Id_Spec}
\end{figure}%
 The input sequence used for identification is analyzed to examine its persistent excitation (PE) properties. 
For example, if we analyze the spectral contents of the input signals, as in Figure \ref{fig:Par_Id_Spec}, we can conclude that, although they are primarily made of low-frequency components, their spectrum is different than zero in a sufficiently wide frequency range. Note that as the sampling time of the datalogger is $\tau_{\rm\scriptscriptstyle id}=10$ s, the maximum frequency in the spectrum is $0.1$ Hz.
\\
An optimization program is defined to minimize the error between the simulated output and the available measurements as follows.
\begin{align*}
\min_{\theta, x(0)} 	&\qquad  \sum_{h=1}^{N} \norm{\tilde{y}(\tau_{\rm\scriptscriptstyle id} h) - \tilde{y}_{\rm\scriptscriptstyle meas}\left[h\right]}\\
\rm{s.t.} 	& \qquad \dot{x}(t) = \tilde{f}(x(t), \tilde{u}_{\rm\scriptscriptstyle meas}\left[\lfloor t/\tau_{\rm\scriptscriptstyle id}\rfloor\right], \theta)\\
& \qquad \tilde{y}(t) = \tilde{g}(x(t), \tilde{u}_{\rm\scriptscriptstyle meas}\left[ \lfloor t/\tau_{\rm\scriptscriptstyle id}\rfloor\right], \theta)\\
& \qquad \theta \in \Theta
\end{align*}
The notation $\lfloor\cdot\rfloor$ denotes the floor function. The decision variable $\theta = [V^{\rm \scriptscriptstyle H}, A^{\rm \scriptscriptstyle B}, A^{\rm \scriptscriptstyle H}, \eta, \beta]$, is a vector  of uncertain model parameters, related to the steam header geometry, the burner and the heat transfer efficiencies. The set $\Theta$ is defined by known physically consistent upper and lower bounds of the parameters to be identified, e.g. while all the unknown parameters are imposed to be positive quantities, efficiencies are also constrained to be strictly less than one. Moreover for each data-set, the initial state is considered as an additional optimization variable.\\
\begin{figure}[bt]
	\centering
	\includegraphics[width=.7\linewidth]{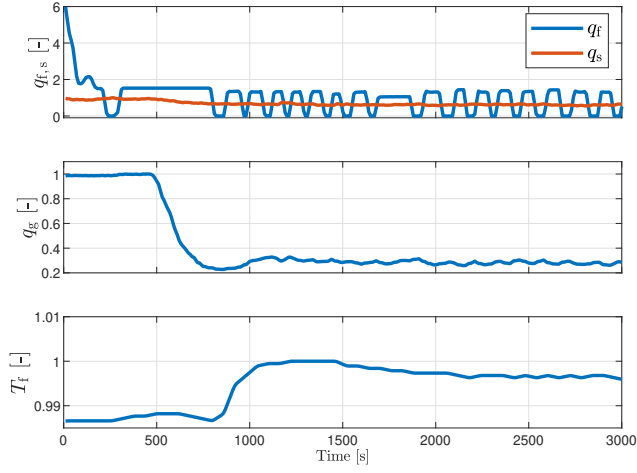}
	\caption{Input data for validation of the identification procedure. Data are adimensionalized for confidentiality}
	\label{fig:Par_Id2}
\end{figure}%
The validation has been performed in a similar process condition, extracted by the available historical data. In  Figure \ref{fig:Par_Id2}, the input data of the validation set are shown. For confidentiality reasons the data shown here are adimensional.
The simulated and measured output trajectories, obtained in the validation phase, are compared in Figure~\ref{fig:Par_Id}.
\subsubsection[]{Thermal stress constraint} \label{subsub:th_model}
One of the main limitations in the boiler operation (and in particular, the start-up) is related to the thermal stress of the shell and the internal tubes. A large thermal stress, due to a too steep increment of the temperatures, leads to a reduction of components' life-cycle, increasing the costs for inspections and maintenance. Therefore, thermal stresses have to be kept under control. The thermal stress $\sigma$ can be modelled as follows.
\begin{equation}
\label{eq:th_mod}
\sigma=k_{\scriptscriptstyle \rm t}(T_{\scriptscriptstyle \rm shell}) \frac{dT^{\scriptscriptstyle \rm shell} }{dt}
\end{equation}
$k_{\scriptscriptstyle \rm t}$ is a function of the material properties and can be either constant or temperature-dependent.\\
For the pressurized components, as the boiler shell, the maximum temperature rate $r_{\scriptscriptstyle \rm T}$ is defined by the European standard EN 12952-3 \cite{EN12952} and must fulfill the following inequality.
\begin{equation}
\label{eq:sig_max}
\abs{\beta_{\scriptscriptstyle \rm p} (p-p_0) \frac{d_{\scriptscriptstyle \rm in} + s}{2s} + \beta_{\scriptscriptstyle \rm T}\frac{E \alpha}{1-\nu} \frac{c_{\scriptscriptstyle \rm p} \rho}{k} s^2 \phi_{\scriptscriptstyle \rm f}  r_{\scriptscriptstyle \rm T}} < \sigma^{\max}
\end{equation}
Variables $p$ and $p_0$ are the nominal and initial, respectively, pressures, $d_{\scriptscriptstyle \rm in}$ is the internal diameter, $s$ the wall thickness, the  mechanical and thermal material properties -$E$, $\nu$, $\alpha$, $\rho$, $k$, and $c_{\scriptscriptstyle \rm p}$- are, respectively, Young's modulus, Poisson's ratio, the thermal expansion coefficient, the density, the thermal conductivity and the specific heat. Also, $\phi_{\scriptscriptstyle \rm f}$ is the cylindrical shape factor, which is a function of the ratio of internal and external diameters.
The first addend is the stress caused by pressurization, while the second one is the thermal stress: $\beta_{\scriptscriptstyle \rm p}$ and $\beta_{\scriptscriptstyle \rm T}$ are the stress concentration factors, respectively for internal pressure and thermal-based stresses. As suggested by~\cite{Taler2015} and references therein, we can set $\beta_{\scriptscriptstyle \rm P}=0.51$ and $\beta_{\scriptscriptstyle \rm T}=2$.\\
The coefficient $k_{\scriptscriptstyle \rm t}$, which can be recovered by inspection from equation~\eqref{eq:sig_max}, is proportional to the Young's modulus and the thermal expansion coefficient, that are in general temperature-dependent. However, for the temperature range considered in this specific application, $E = 1.82\, 10^5 $MPa and $\alpha = 1.35\, 10^{-5} \text{m}^2/\text{s}$ and $k_{\scriptscriptstyle \rm t}$ can be considered constant.\\
The thermal stress inequality~\eqref{eq:sig_max} can be formulated as a constraint on the rate of change of the component temperature. In particular, in the model predictive control optimal control problem, we will enforce the following inequality on the difference between the wall temperature values at subsequent sampling times:
\begin{equation}
\label{eq:th_bound}
(T_{\rm \scriptscriptstyle w}[h+1]-T_{\rm \scriptscriptstyle w}{\scriptscriptstyle j}[h])/\tau \leq  r_{\scriptscriptstyle \rm{T}}^{\max}
\end{equation}
where $\tau$ denotes the used sampling time. The value of $r_{\scriptscriptstyle \rm T}^{\max}$ is computed from equation~\eqref{eq:sig_max} by imposing the maximum allowable stress. The bound~\eqref{eq:th_bound} is imposed to $T_{\rm\scriptscriptstyle w}$ under the assumption that the latter is equal to the water temperature, which in turn is assumed to be in equilibrium with the shell one. An analogous constraint on temperature variation must be imposed on internal tubes, where we express the \eqref{eq:th_bound} inequality considering respectively $T^{\rm\scriptscriptstyle t,GAS}$ and $T^{\rm\scriptscriptstyle t,CHP}$ for the main burner tubes and the one connected to CHP exhaust. 
\subsection{The CHP high-level hybrid model}
\label{sec:Model-sub:CHP_H}
\begin{figure}[hbt]
	\centering

	\def\svgwidth{0.5\columnwidth}

	\includegraphics[width=.5\linewidth]{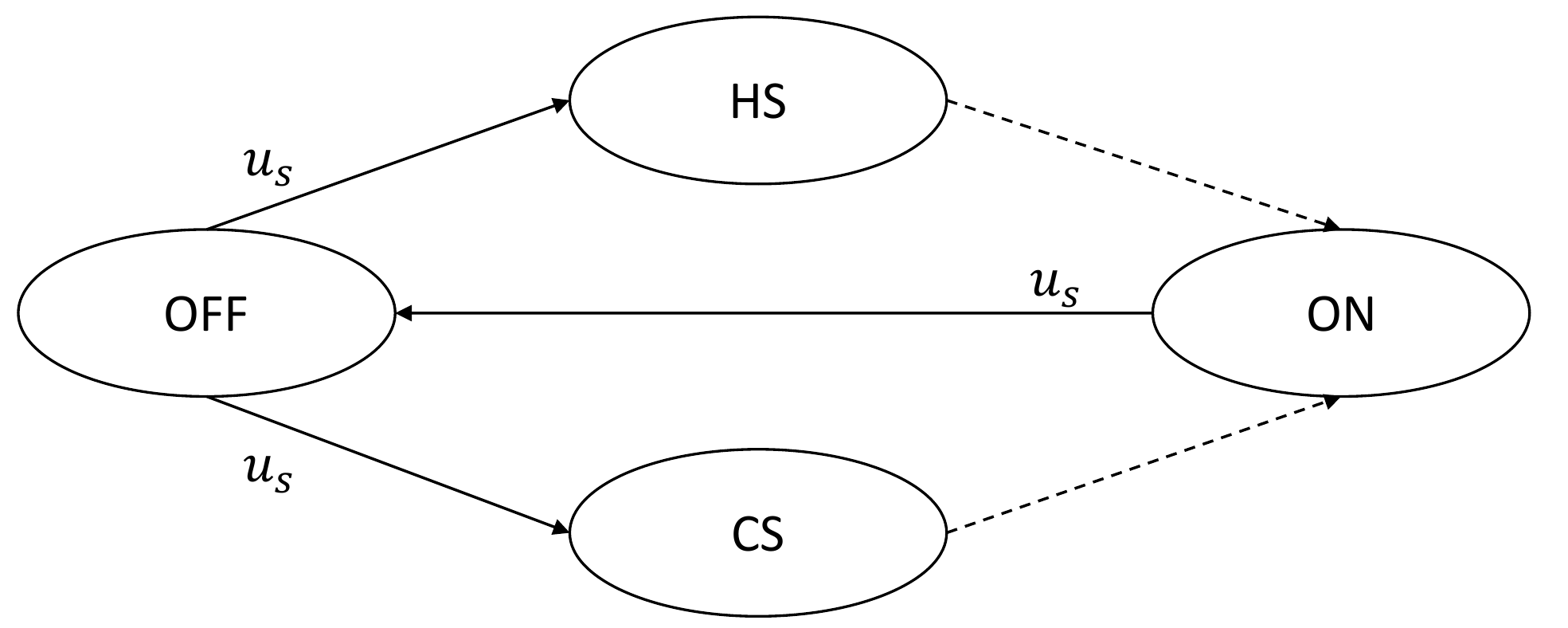}
	\caption{Finite State Machine of the internal combustion engine. Di-graph represents the system modes and the possible transitions. Only decision variables are shown, however transitions are also implied by dwell time conditions. Solid edges denote externally governed transitions, while dotted lines show transitions induced by state guard conditions}

	\label{fig:States}
\end{figure}

The model used by the high-level optimizer, having a sampling time $\tau_{\rm\scriptscriptstyle H}$, must describe the steady-state behaviour of the CHP in all its operation modes and the transitions among them. To this end, a discrete hybrid model is used, which describes the time-evolution of a set of continuous and discrete states, subject to both discrete (boolean) and continuous inputs.\\
The discrete states, $m^{\rm \scriptscriptstyle CHP}$, correspond to the different CHP operating modes - shut down (OFF), cold start (CS), hot start (HS), and production (ON) - thus we define the discrete state to lie in the set $\mathcal{M}$: $m^{\rm \scriptscriptstyle CHP} \in \{\rm{OFF, CS, HS, ON}\}$.\\ The high-level model is based on the assumption that the evolution among the operating modes can be described by temporal thresholds.\\ 
The dynamic state is  therefore an interal clock, denoted with $\tau^{\rm\scriptscriptstyle CHP}$, that represents the number of sampling times spent by the system in the present operation mode, where $\tau^{\rm\scriptscriptstyle CHP}\in \mathbb{Z}_0^+ $.\\
A specific class of discrete hybrid automata (DHA) is the timed automata (TA) \cite{Lygeros2008}, where the discrete states and transitions are governed by  dwell-times and defined by timed transition languages.\\
The continuous input, based on the closed-loop model described above, is the demanded (and produced) electric power $P_{\rm \scriptscriptstyle e}^{\rm \scriptscriptstyle CHP}$, while the model output is the pair $(q_{\rm \scriptscriptstyle g}^{\rm \scriptscriptstyle CHP},H_{\rm\scriptscriptstyle ex}^{\rm\scriptscriptstyle CHP})$ which depends on the discrete state $m^{\rm\scriptscriptstyle CHP}$ and on the continuous input $P_{\rm \scriptscriptstyle e}^{\rm \scriptscriptstyle CHP}$.\\
In particular, in the operating mode, i.e. $m^{\rm \scriptscriptstyle CHP}=\rm{ON}$, the output pair of the CHP is given by affine relationships with input $P_{\rm \scriptscriptstyle e}^{\rm \scriptscriptstyle CHP}$:
\begin{subequations}\label{eq:CHP_outputs_HL_tot}
\begin{align}
\begin{array}{lcl}
q_{\rm \scriptscriptstyle g}^{\rm \scriptscriptstyle CHP}&=&(P_{\rm \scriptscriptstyle e}^{\rm \scriptscriptstyle CHP}- P_{\rm\scriptscriptstyle int})/\gamma_{q}\\
H_{\rm\scriptscriptstyle ex}^{\rm\scriptscriptstyle CHP}&=&\gamma_{h}P_{\rm \scriptscriptstyle e}^{\rm \scriptscriptstyle CHP}
\end{array}\label{eq:IN-OUT-CHP}
\end{align}
$P_{\rm\scriptscriptstyle int}$, $\gamma_{q}$, and $\gamma_h$ are identified from data, as shown in Figure~\ref{fig:Static_CHP}.\\
For what concern all the other operating modes, the outputs of the CHP are both null, while a constant gas consumption is considered in both starting modes,  $\bar{q}_{\rm \scriptscriptstyle g \, CS}^{\rm \scriptscriptstyle CHP}$ and $\bar{q}_{\rm \scriptscriptstyle g \, HS}^{\rm \scriptscriptstyle CHP}$ respectively for the cold start and hot start mode.
\begin{align}(q_{\rm \scriptscriptstyle g}^{\rm \scriptscriptstyle CHP},H_{\rm\scriptscriptstyle ex}^{\rm\scriptscriptstyle CHP})=
\left\{\begin{array}{ll}
(0,0)&\text{if }m^{\rm \scriptscriptstyle CHP}=\rm{OFF}\\
(\bar{q}_{\rm \scriptscriptstyle g \, CS}^{\rm \scriptscriptstyle CHP},0)&\text{if }m^{\rm \scriptscriptstyle CHP}=\rm{CS}\\
(\bar{q}_{\rm \scriptscriptstyle g \, HS}^{\rm \scriptscriptstyle CHP},0)&\text{if }m^{\rm \scriptscriptstyle CHP}=\rm{HS}
\end{array}
\right.\label{eq:CHP_outputs_HL}\end{align}
\end{subequations}
Regarding the transition between the operating modes, the feasible changeovers are described by the directed graph depicted in Figure~\ref{fig:States}. As introduced, these transitions are possible only after a defined dwell time, i.e. if the time spent in the current operating mode $\tau^{\rm\scriptscriptstyle CHP}(k)$ is greater than a proper threshold.
In addition, the $\rm{ON}\rightarrow \rm{OFF}$ transition and the ones involving the system startup ($\rm{OFF}\rightarrow \rm{CS}$ and $\rm{OFF}\rightarrow \rm{HS}$) are triggered when a proper switching input signal $u_{\rm\scriptscriptstyle S}^{\rm\scriptscriptstyle CHP}(k)\in\{0,1\}$ is set to $1$. \\

The high-level model of the system is therefore expressed in the DHA framework. An open DHA model is described by the tuple  
$ \mathcal{H} = (\mathcal{M},\mathcal{X}, \mathcal{U}, \mathcal{W}, f, h, \mathit{Init},\mathit{Dom},\mathcal{E},\mathcal{G},\mathcal{R})$. The discrete modes are defined by the finite set $\mathcal{M}$, while the continuous variables are given by the finite set of real-valued variables $\mathcal{X}$, as in \eqref{eq:DHA_qxsets}. The sets $\mathcal{U}$ and $\mathcal{W}$ include the input and output variables, see \eqref{eq:DHA_uysets}, considering both continuous and discrete ones. The vector field $f$ expresses the evolution of continuous variables and it generally depends on the active discrete mode $m$, as for the output map $h$. The affine equations in \eqref{eq:DHA_h} depend indeed on the active mode, and are detailed in \eqref{eq:CHP_outputs_HL_tot}.
While $\mathit{Init}$ is the set of valid initial conditions, $\mathit{Dom}$ represents the domain in which the discrete mode $m$ is invariant. The possible transitions between these modes are defined by the edges $\mathcal{E}$, and transitions are activated if the guard conditions - described by $\mathcal{G}$ -  are met. The transition might involve a reset of the continuous variable, which is defined in \eqref{eq:DHA_Reset} by $\mathcal{R}$. Here in detail, the model of the CHP in the DHA framework is reported:\\
\begin{subequations}
	\begin{align}
	\mathcal{M} = \{\rm{OFF, CS, HS, ON} \} \quad&\quad  \mathcal{X} = \{\tau^{\rm\scriptscriptstyle CHP}\} \in \mathbb{Z}_0^+
	\label{eq:DHA_qxsets}
	\end{align}
		\begin{align}
	\mathcal{U}= \{ P_{\rm \scriptscriptstyle El}^{\rm \scriptscriptstyle CHP},  u_{\rm\scriptscriptstyle S}^{\rm\scriptscriptstyle CHP}  \} \in \mathbb{R} \times \{0,1\} \quad&\quad \mathcal{W}=\mathbb{R}^2 
	\label{eq:DHA_uysets}
	\end{align}
	\begin{align}
		f: &\mathcal{M}\times \mathcal{X} \times \mathcal{U} \to \mathcal{X}  &\quad f(m^{\rm\scriptscriptstyle CHP},  \tau^{\rm\scriptscriptstyle CHP}, u^{\rm\scriptscriptstyle CHP})=  \tau^{\rm\scriptscriptstyle CHP} + 1 &\quad \forall m^{\rm\scriptscriptstyle CHP} \in \mathcal{M} \label{eq:DHA_f}\\
	h: &\mathcal{M}\times \mathcal{X} \times \mathcal{U} \to \mathcal{W}  &\quad 	 h(m^{\rm\scriptscriptstyle CHP},  \tau^{\rm\scriptscriptstyle CHP}, u^{\rm\scriptscriptstyle CHP}) = \text{aff}(u^{\rm\scriptscriptstyle CHP})& \label{eq:DHA_h}
	\end{align}
	\begin{equation}
	\mathcal{E} = \{\rm{(OFF, CS),  (OFF,HS), (CS,ON), (HS,ON), (ON,OFF)}\}
	\end{equation}
\begin{align}
\mathit{Dom}(m) = \left\{\begin{array}{l}
\{ \tau^{\rm\scriptscriptstyle CHP}\leq \tau_{\scriptscriptstyle \rm{OFF}\rightarrow \rm{HS}}\}\vee \{u_{\rm\scriptscriptstyle S}^{\rm\scriptscriptstyle CHP}(k)=0\}\\
\{ \tau^{\rm\scriptscriptstyle CHP}\leq \tau_{\scriptscriptstyle \rm{CS}\rightarrow \rm{ON}}\}\\
\{ \tau^{\rm\scriptscriptstyle CHP}\leq \tau_{\scriptscriptstyle \rm{HS}\rightarrow \rm{ON}}\}\\
\{ \tau^{\rm\scriptscriptstyle CHP}\leq \tau_{\scriptscriptstyle \rm{ON}\rightarrow \rm{OFF}}\}\vee \{u_{\rm\scriptscriptstyle S}^{\rm\scriptscriptstyle CHP}(k)=0\}
\end{array}
\right. \quad& \begin{array}{l}
\text{if } m^{\rm\scriptscriptstyle CHP} = \rm{OFF} \\
\text{if } m^{\rm\scriptscriptstyle CHP} = \rm{CS} \\
\text{if } m^{\rm\scriptscriptstyle CHP} = \rm{HS} \\
\text{if } m^{\rm\scriptscriptstyle CHP} = \rm{ON} 
\end{array}
\label{eq:DHA_Domset}
\end{align}
\begin{align}
\mathcal{G}(e) = \left\{\begin{array}{l}
\{\tau^{\rm\scriptscriptstyle CHP}\in (\tau_{\scriptscriptstyle \rm{OFF}\rightarrow \rm{HS}},\tau_{\scriptscriptstyle \rm{OFF}\rightarrow \rm{CS}}]\}\wedge \{u_{\rm\scriptscriptstyle S}^{\rm\scriptscriptstyle CHP}(k)=1\}\\
\{ \tau^{\rm\scriptscriptstyle CHP}>\tau_{\scriptscriptstyle \rm{OFF}\rightarrow \rm{CS}}\}\wedge \{u_{\rm\scriptscriptstyle S}^{\rm\scriptscriptstyle CHP}(k)=1\}\\
\{ \tau^{\rm\scriptscriptstyle CHP}>\tau_{\scriptscriptstyle m\rightarrow \rm{ON}}\}\\
\{ \tau^{\rm\scriptscriptstyle CHP}>\tau_{\scriptscriptstyle \rm{ON}\rightarrow \rm{OFF}}\}\wedge \{u_{\rm\scriptscriptstyle S}^{\rm\scriptscriptstyle CHP}(k)=1\}
\end{array}
\right. \quad& \begin{array}{l}
\text{if } e^{\rm\scriptscriptstyle CHP} = (\rm{OFF, HS}) \\
\text{if } e^{\rm\scriptscriptstyle CHP} = (\rm{OFF, CS}) \\
\text{if } e^{\rm\scriptscriptstyle CHP} = (\rm{CS}, \rm{ON})\vee(\rm{HS}, \rm{ON}) \\
\text{if } e^{\rm\scriptscriptstyle CHP} = (\rm{ON, OFF})
\end{array}
 \label{eq:DHA_Guard}
\end{align}
\begin{equation}\label{eq:DHA_Reset}
\mathcal{R}(e^{\rm\scriptscriptstyle CHP},  \tau^{\rm\scriptscriptstyle CHP}, u^{\rm\scriptscriptstyle CHP}) = 0  \quad\quad \forall e^{\rm\scriptscriptstyle CHP} \in \mathcal{E}
\end{equation}
\label{eq:DHA}
\end{subequations}

%
The values of $\tau_{\scriptscriptstyle \rm{OFF}\rightarrow \rm{HS}}$, $\tau_{\scriptscriptstyle \rm{OFF}\rightarrow \rm{CS}}$, $\tau_{\scriptscriptstyle \rm{CS}\rightarrow \rm{ON}}$, $\tau_{\scriptscriptstyle \rm{HS}\rightarrow \rm{ON}}$, and $\tau_{\scriptscriptstyle \rm{ON}\rightarrow \rm{OFF}}$, are suitably-defined thresholds, which characterize the time abstraction considered at the high level for the operating modes.
To delve deeper into the essence of these thresholds, a distinction has to be done between them:
$\tau_{\scriptscriptstyle \rm{ON}\rightarrow \rm{OFF}}$ is a value imposed by the user to avoid an immediate shutdown after the system activation;
$\tau_{\scriptscriptstyle \rm{OFF}\rightarrow \rm{HS}}$ is determined by the minimum time requirements for the shutdown process; the other thresholds are instead related to the modeling approximation.\\
 As a matter of fact, the duration of the startup time is not a fixed value, but is related to the dynamic behavior of continuous variables and most of all depends on their initial values. However, as a modeling design choice, the startup duration has been clustered in two possible classes, cold startup and hot startup, which are characterized by two durations, i.e., respectively $\tau_{\scriptscriptstyle \rm{CS}\rightarrow \rm{ON}}$ and $\tau_{\scriptscriptstyle \rm{HS}\rightarrow \rm{ON}}$, due to the low/high value of the initial state. This status is conditioned by the time spent in the OFF mode.

\subsection{The boiler high-level hybrid model}
\label{sec:Model-sub:boiler_H}
\begin{figure}[thpb]
	\centering
	\includegraphics[width=.7\linewidth]{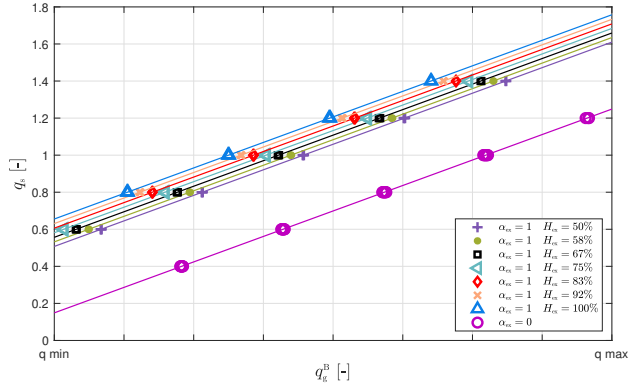}
	\caption{Relationship between $q_{\rm \scriptscriptstyle g}^{\rm \scriptscriptstyle B}$ and $q_{\rm \scriptscriptstyle s}^{\rm \scriptscriptstyle B}$ in steady-state conditions. The identified static map \eqref{eq:Boiler_Stat} is presented for different values of $, u_{\rm\scriptscriptstyle ex} H_{\rm \scriptscriptstyle ex}$ in solid lines. Markers show the steady-state data of the nonlinear model \eqref{eqn:cons_eq_in_bis}.}
	\label{fig:Static_Boiler}
\end{figure}

The boiler model for the high level optimizer is also devised as a hybrid model which describing the steady-state operation of the FTB for each functional mode.\\
\begin{figure}[thpb]
	\centering
	\includegraphics[width=.5\linewidth]{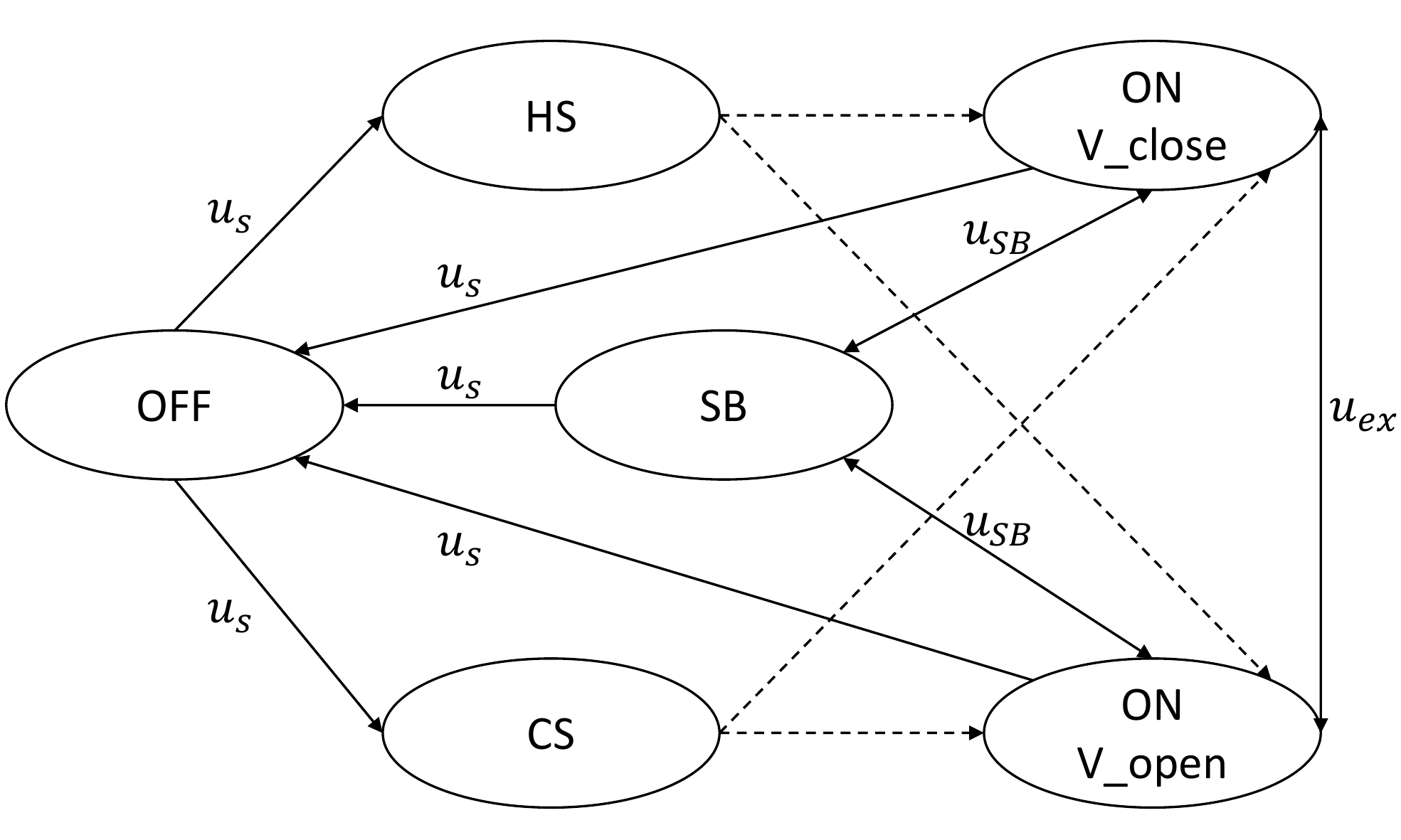}
	\caption{Finite State Machine of the fire tube boiler. Di-graph represents the system modes and the possible transitions. Only decision variables are shown, however transitions are also implied  by dwell time conditions}
	\label{fig:DHA_Boiler}
\end{figure}%
The Finite State Machine of the FTB, see Figure \ref{fig:DHA_Boiler}, extends the one just presented for the CHP, as the boiler has two production states - because of the interlink between CHP and FTB - and a Stand-by mode. The operating production mode can be actively selected by commanding the exhaust valve, $u_{\rm\scriptscriptstyle ex}$, which can switch between two discrete values (open/closed). When open, the CHP flue gases are diverted inside the boiler.\\
Moreover, as discussed in \ref{sec:stand_by}, the FTB can operate in any of the two productive modes only between a minimum/maximum steam generation, where $\bar{q}_{\rm \scriptscriptstyle s}^ {\scriptscriptstyle \min} > 0$.
To address this condition, the FTB includes another additional mode, which is a Stand-By (SB) status commanded by the boolean input $u_{\rm\scriptscriptstyle SB}^{\rm\scriptscriptstyle B}$. 
In this manner, the discrete state $m^{\rm\scriptscriptstyle B}$ can lie in the set $\mathcal{M}^{\scriptscriptstyle B}= \{\rm{OFF, CS, HS, SB}, \rm{ON}_{u_{\rm\scriptscriptstyle ex}=0},\rm{ON}_{u_{\rm\scriptscriptstyle ex}=1}\}$.\\
The high level model therefore presents two possible alternative  reachable states when the steam demand is not required: the OFF and the SB modes. The choice between the two determines a different evolution of the FSM to return on the operative modes and it has a dissimilar impact on system consumption.\\
The SB mode, detailly discussed in Section~\ref{sec:stand_by}, has a low but non-zero consumption of natural gas, while it presents the advantage of being able to return to productive mode almost instantaneously, just activating the controllable variable  $u_{\rm\scriptscriptstyle SB}^{\rm\scriptscriptstyle B}$.
Instead, switching off the boiler when the consumer demand is zero, determines a null fuel consumption during idle. As a drawback, the FTB needs to execute the start-up procedures before being ready to provide steam: thus the FTB cannot be responsive to a sudden steam demand. The activation command, $u_{\rm\scriptscriptstyle s}$ must consider the mandatory dwell time related to the CS (or  HS) states.\\
The model, in the productive modes $\rm{ON}_{u_{\rm\scriptscriptstyle ex}=0}$ and $\rm{ON}_{u_{\rm\scriptscriptstyle ex}=1}$, is assumed to be in steady-state conditions and under control, in such a way that $p_{\scriptscriptstyle \rm s}=\bar{p}_{\scriptscriptstyle \rm s}$ and $l_{\scriptscriptstyle \rm w}=\bar{l}_{\scriptscriptstyle \rm w}$, where $\bar{p}_{\scriptscriptstyle \rm s}$ and $\bar{l}_{\scriptscriptstyle \rm w}$ are suitable and fixed nominal conditions. In this way, the combustion gas flow-rate required to guarantee the steam demand $q^{\scriptscriptstyle \rm B}_{\rm\scriptscriptstyle g}$ is a function of $q_{\scriptscriptstyle \rm s}^{\scriptscriptstyle \rm B}$ and, when $u_{\rm\scriptscriptstyle ex}=1$, of the flue gas enthaply too, $H^{\scriptscriptstyle \rm CHP}_{\scriptscriptstyle \rm ex}$. For simplicity, an affine relationship has been identified and considered in the higher control layer:
\begin{align} \label{eq:Boiler_Stat}
q^{\scriptscriptstyle \rm B}_{\rm\scriptscriptstyle g}=\gamma_1 q_{\scriptscriptstyle \rm s}^{\scriptscriptstyle \rm B}+\gamma_2 u_{\rm\scriptscriptstyle ex}H^{\scriptscriptstyle \rm CHP}_{\scriptscriptstyle \rm ex}+\gamma_3
\end{align}
where $\gamma_1$, $\gamma_2$, and $\gamma_3$ are identified from the model static map, as shown in Figure \ref{fig:Static_Boiler}.\\
In all the other operating modes, the produced steam flow rate is null or below $q_{\rm \scriptscriptstyle s}^ {\scriptscriptstyle \min}$, while a natural gas consumption is present in SB, CS and HS modes:
 \begin{align}(q_{\rm \scriptscriptstyle g}^{\rm \scriptscriptstyle B},q_{\rm \scriptscriptstyle s}^{\rm \scriptscriptstyle B})=
 \left\{\begin{array}{ll}
 (0,0)&\text{if }m^{\rm \scriptscriptstyle B}=\rm{OFF}\\
 (\bar{q}_{\rm \scriptscriptstyle g \, CS}^{\rm \scriptscriptstyle B},0)&\text{if }m^{\rm \scriptscriptstyle B}=\rm{CS}\\
 (\bar{q}_{\rm \scriptscriptstyle g \, HS}^{\rm \scriptscriptstyle B},0)&\text{if }m^{\rm \scriptscriptstyle B}=\rm{HS}\\
 (\bar{q}_{\rm \scriptscriptstyle g \, SB}^{\rm \scriptscriptstyle B},\bar{q}_{\rm \scriptscriptstyle s \, SB}^{\rm \scriptscriptstyle B})&\text{if }m^{\rm \scriptscriptstyle B}=\rm{SB}
 \end{array}
 \right.\label{eq:SB_outputs_HL}\end{align}
 where $\bar{q}_{\rm \scriptscriptstyle g \,CS}^{\rm \scriptscriptstyle B}$ and $\bar{q}_{\rm \scriptscriptstyle g \, HS}^{\rm \scriptscriptstyle B}$ are given values of gas consumption for the start-up modes, while $\bar{q}_{\rm \scriptscriptstyle g \, SB}^{\rm \scriptscriptstyle B},\bar{q}_{\rm \scriptscriptstyle s \, SB}^{\rm \scriptscriptstyle B}$ are the average gas and steam flow-rates of the limit cycle occurring in the stand-by mode.
Transitions are possible only if the time spent in the current operating mode $\tau^{\rm\scriptscriptstyle B}$ is greater than a proper threshold: in particular the dwell-time of the stand-by mode is null, $\tau^{\rm\scriptscriptstyle B}_{\rm\scriptscriptstyle SB}=0$, therefore FTB can be very reactive to improvise steam demand, while in stand-by condition.\\
Analogously to what discussed above for the high-level model of the CHP, the timed transitions governing hot/cold startup are exclusively a modeling construction. The actual duration of the startup depends on the initial condition of the FTB, in particular the water temperature, while its rate of change is limited by the thermal stress constraints, as discussed in Section \ref{chap:ltv_mpc}. The quick reactivity of the stand-by mode is substantially due to the fact that this operating mode maintains the water temperature very close the its set-point. Instead, when the FTB is shut down, the temperature decreases with an evolution well described by a first order dynamics. The internal temperature and therefore the initial condition of the start-up procedure are almost entirely related on the time spent by the FTB in OFF mode. Consequently, the duration of the start-up will be longer, the longer the system has been shut off.
To simplify the model, two possible start-up times have been defined, clustering the behaviour of the system just in cold start-up and hot start-up: 
$\tau_{\scriptscriptstyle \rm{OFF}\rightarrow \rm{CS}}$ discriminates if the unit has been in OFF mode enough to experience a water temperate below the chosen threshold, while the dwell times for cold/hot startup modes are over-approximations of the actual cold and hot start-up times.\\
Also, both the switch-off and switch-on transitions are triggered when a proper switching input signal $u_{\rm\scriptscriptstyle S}^{\rm\scriptscriptstyle B}[k]\in\{0,1\}$ is set to $1$. As said,  a binary command $u_{\rm\scriptscriptstyle SB}^{\rm\scriptscriptstyle B}[k]\in\{0,1\}$ is introduced to move the boiler in and out the stand-by mode.\\
Therefore, the input vector of the high level model is $u^{\rm\scriptscriptstyle B} = \left\lbrace  q_{\rm \scriptscriptstyle s}^{\rm \scriptscriptstyle B}, H^{\scriptscriptstyle \rm CHP}_{\scriptscriptstyle \rm ex}, u_{\rm\scriptscriptstyle S}^{\rm\scriptscriptstyle B}, u_{\rm\scriptscriptstyle SB}^{\rm\scriptscriptstyle B}, u_{\rm\scriptscriptstyle ex} \right\rbrace $, where $u \in \mathcal{U} \subseteq \mathbb{R}^2_{\scriptscriptstyle \geq 0} \times \{0,1\}^3$.

%

\section{Hierarchical control}
\label{chap:hier_ctrl}
\begin{figure}[!ht]
	\centering
	\includegraphics[width=1\linewidth]{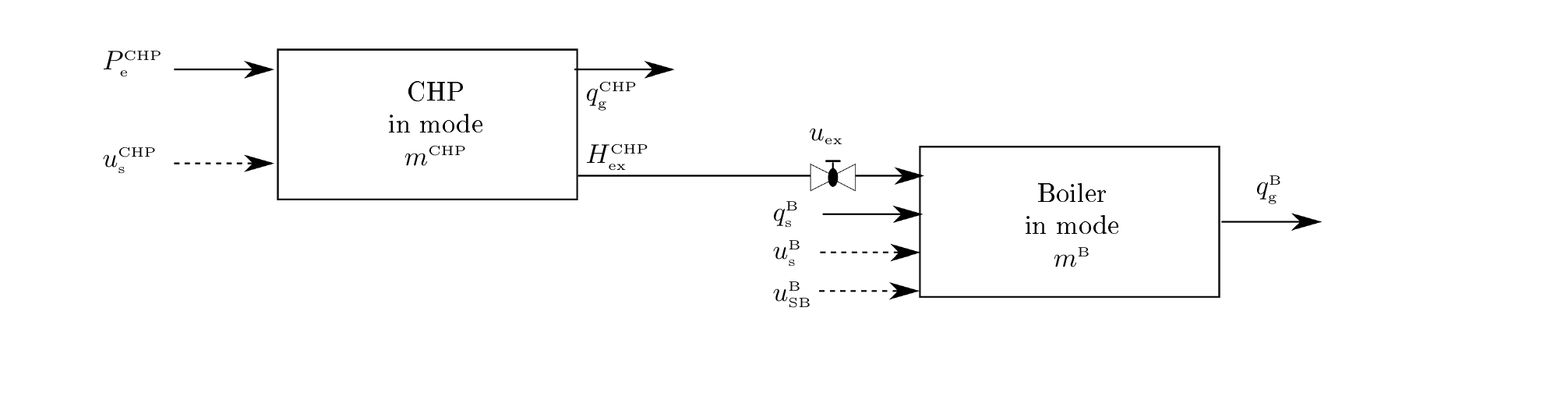}
	\caption{Sketch of the high-level model}
	\label{fig:HLmodel}
\end{figure}

In this section we propose a control scheme that allows to satisfy the electricity and heat demands in an economically optimal way, at minimum operative cost, considering the combined heat and power generation unit composed of the CHP and the boiler described in Section~\ref{chap:model}.
The demand of steam and electricity is considered given and known for a $24$-hours horizon, as well as a price forecast. We enable the possibility to exchange power with the electrical grid: the surplus of electricity produced by the CHP can be sold, or, if negative, it can be bought from the market.\\
The hierarchical control structure here proposed includes a high optimization layer that aims to optimize the performance of the integrated plant on a day-ahead basis with a sampling time of $\tau_{\rm\scriptscriptstyle H}=15$ min. The model used at the high level is sketched in Figure~\ref{fig:HLmodel}.
The high level optimizer receives as input the 24-hours ahead electric and steam demands $P_{\rm\scriptscriptstyle e}^{\rm\scriptscriptstyle  D}[k]$ and $q_{\rm\scriptscriptstyle s}^{\rm\scriptscriptstyle D}[k]$, respectively, for $k=1,\dots, 24\times 4$, together with a forecast of the electricity prices. Based on this, the optimizer determines the future modes of operation of the systems $m^{\rm\scriptscriptstyle CHP}[k]$ and $m^{\rm\scriptscriptstyle B}[k]$, as well as the optimal production profiles $P_{\rm\scriptscriptstyle e}^{\rm\scriptscriptstyle  CHP}[k]$, $q_{\rm\scriptscriptstyle s}^{\rm\scriptscriptstyle B}[k]$ and the corresponding predicted combustion gas amounts $q_{\rm\scriptscriptstyle g}^{\rm\scriptscriptstyle CHP}[k]$ and $q_{\rm\scriptscriptstyle g}^{\rm\scriptscriptstyle B}[k]$ for the whole future optimization time horizon $k=1,\dots,24\times 4$.\\
The lower level is a fast regulating one with a sampling time of $\tau_{\rm\scriptscriptstyle L}=60$ s, aiming to control the dynamic variables of the subsystems.\\
For clarity, the symbol $h$ will be used to denote the discrete time index in the fast time scale of the lower layer, while $k$ will be used to denote the time index in the slow high-level time scale.

\subsection{The high level optimization problem}
\label{sec:control-sub:HL}
The high level optimizer objective is the minimization of the operating costs of the plant on a windows of 24-hours, while  fulfilling completely the electric and steam demands in the context of a liberalized energy market.
In addition to the subsystem operating points, a further decision variable is considered in the high-level optimization, by including the possibility either to buy electricity $P_{\rm \scriptscriptstyle e}^{\rm \scriptscriptstyle Purch }$ from the power grid, to support the production of the CHP unit or  to sell surplus of produced energy to the grid.\\
The need to consider, in an integrated fashion, both the ICE and the FTB originates from the fact that a valve is present to divert exhaust gases from the CHP to a set of tubes in the boiler. Note that it can be opened only when the engine is in production, which translate in the following logical implication:
\begin{equation}
\label{eqn:log_cond}
u_{\rm \scriptscriptstyle ex} = 1 \implies m^{\rm \scriptscriptstyle CHP} = \rm{ON}
\end{equation}
where $ m^{\rm \scriptscriptstyle CHP}$ is the operating mode of the internal combustion engine.
To formulate the high-level optimization problem, the two hybrid automata devised in Sections~\ref{sec:Model-sub:CHP_H} and~\ref{sec:Model-sub:boiler_H} have been integrated in a single Mixed Logical Dynamical (MLD) system~\cite{Bemporad1999} of the type:
\begin{align}
&x_{\rm\scriptscriptstyle H}[k+1] =  A_{\rm\scriptscriptstyle H} x_{\rm\scriptscriptstyle H}[k] + B_{\rm\scriptscriptstyle H,1} u_{\rm\scriptscriptstyle H}[k] + B_{\rm\scriptscriptstyle H,2} \Delta_{\rm\scriptscriptstyle H}[k] + B_{\rm\scriptscriptstyle H,3} z_{\rm\scriptscriptstyle H}[k] \nonumber\\
&y_{\rm\scriptscriptstyle H}[k] = C_{\rm\scriptscriptstyle H} x_{\rm\scriptscriptstyle H}[k] + D_{\rm\scriptscriptstyle H,1} u_{\rm\scriptscriptstyle H}[k] + D_{\rm\scriptscriptstyle H,2} \Delta_{\rm\scriptscriptstyle H}[k] + D_{\rm\scriptscriptstyle H,3} z_{\rm\scriptscriptstyle H}[k]
\label{eq:MLD}\\
&E_{\rm\scriptscriptstyle H,4} \Delta_{\rm\scriptscriptstyle H}[k] + E_{\rm\scriptscriptstyle H,5} z_{\rm\scriptscriptstyle H}[k] \leq  E_{\rm\scriptscriptstyle H,1} u_{\rm\scriptscriptstyle H}[k] + E_{\rm\scriptscriptstyle H,4} x_{\rm\scriptscriptstyle H}[k] + E_{\rm\scriptscriptstyle H,5} \nonumber
\end{align}
where the state, ${x}_{\rm\scriptscriptstyle H} = [\tau^{\rm \scriptscriptstyle CHP  },  {\chi}^{\rm \scriptscriptstyle CHP  }_b, \tau^{\rm \scriptscriptstyle B  }, {\chi}^{\rm \scriptscriptstyle B  }_b]'$, is composed of real and boolean variables. More specifically, ${\chi}_b$ is a boolean vector of suitable size to translate the discrete modes $m^{\rm \scriptscriptstyle CHP  }$ and $m^{\rm \scriptscriptstyle B}$ in a binary representation.  The input of the MLD system $u_{\rm\scriptscriptstyle H}$$=[u_{\rm \scriptscriptstyle S}^{\rm \scriptscriptstyle CHP  },$ $ P_{\rm \scriptscriptstyle e}^{\rm \scriptscriptstyle CHP  }, $ $ u_{\rm \scriptscriptstyle S}^{\rm \scriptscriptstyle B  },$ $ u_{\rm \scriptscriptstyle SB}^{\rm \scriptscriptstyle B  },$ $ q_{\rm \scriptscriptstyle s}^{\rm \scriptscriptstyle B  }, $ $u_{\rm \scriptscriptstyle ex},$ $ P_{\rm \scriptscriptstyle e}^{\rm \scriptscriptstyle  Purch }]$ is composed of, respectively, the switching signals and the continuous inputs of the ICE and FTB subsystems, the divert valve command and the electricity purchased from the grid. The output  is $y_{\rm\scriptscriptstyle H} =  [P_{\rm \scriptscriptstyle e}, q^{\rm \scriptscriptstyle CHP}_{\rm \scriptscriptstyle g}, q_{\rm \scriptscriptstyle g}^{\rm\scriptscriptstyle B}]$, where $P_{\rm \scriptscriptstyle e} = P_{\rm \scriptscriptstyle e}^{\rm \scriptscriptstyle CHP} + P_{\rm \scriptscriptstyle e}^{\rm \scriptscriptstyle  Purch }$ is the electric power both purchased and generated by the ICE and the other are the gas consumed  by the systems. $\Delta_{\rm\scriptscriptstyle H}$ and $z_{\rm\scriptscriptstyle H}$ are, respectively, auxiliary binary and continuous variables.
Note that the last equation in \eqref{eq:MLD} includes both the inequalities related to the logical relationships that regulate the switches between discrete states, and the operational constraints \eqref{eq:CHP_CONSTRAINT} and \eqref{eq:boiler_CONSTRAINT}.
The cost function is given by:
\begin{equation}
J_{\rm\scriptscriptstyle H} = C_{\rm g} - R
\label{eq:cost_HL}
\end{equation}
where the fuel cost $C_{\rm g}$ is related to fuel consumption used to respect the utility demand as follows:
\begin{equation}
C_{\rm g} = \sum_{k} (q_{\rm \scriptscriptstyle g}^{\rm \scriptscriptstyle CHP  }[k] + q_{\rm \scriptscriptstyle g}^{\rm \scriptscriptstyle B  }[k]) \lambda_{\rm g}
\label{eq:cost_HL_g}
\end{equation}
The cost~\eqref{eq:cost_HL_g} considers also the costs related to the start-up modes $\{\rm CS, HS\}$ and stand-by, due to the gas flow values defined in \eqref{eq:CHP_outputs_HL} and \eqref{eq:SB_outputs_HL} also for the non productive phases.
The revenue term $R$ in~\eqref{eq:cost_HL} is related to the possibility of selling the electricity surplus to the grid. It can  take into account  also the purchase of electricity, in case $R <0$. By splitting the deviation from the electrical demand in positive and negative fluctuation, $\Delta^\uparrow e$, $\Delta^\downarrow e$, it is possible to consider different prices for electricity, i.e., whether it is sold or bought:
\begin{align}
&R = \sum_{k} \Delta^\uparrow e[k]  \lambda_{\uparrow} [k] - \Delta^\downarrow e[k] \lambda_{\downarrow} [k]\\
\rm{with}& \nonumber\\
&\Delta e[k] = P_{\rm \scriptscriptstyle e}^{\rm \scriptscriptstyle CHP  } [k] +  P_{\rm \scriptscriptstyle e}^{\rm \scriptscriptstyle Purch  } [k] - P_{\rm \scriptscriptstyle e}^{\rm \scriptscriptstyle D  } [k]\nonumber\\
&\Delta e[k] =  \Delta^\uparrow e[k]  - \Delta^\downarrow e[k]\nonumber\\
&\Delta^\uparrow e[k] , \Delta^\downarrow e[k] \geq 0 \quad \forall k \nonumber
\end{align}
The optimization problem can be stated as follows, where the initial time is $k=1$ for simplicity reasons:
\begin{align*}
&\min_{\xi(1) \dots \xi(N_{\rm\scriptscriptstyle H})}  \quad J_{\rm\scriptscriptstyle H} \\
\text{s.t.} & \quad \text{MLD model~\eqref{eq:MLD}}\\
& \quad P_{\rm \scriptscriptstyle e}^{\rm \scriptscriptstyle CHP  } [k] + P_{\rm \scriptscriptstyle e}^{\rm \scriptscriptstyle Purch  } [k]  \geq P_{\rm \scriptscriptstyle e}^{\rm \scriptscriptstyle D  } [k]\\
& \quad q_{\rm \scriptscriptstyle S}^{\rm\scriptscriptstyle B}[k] \geq q_{\rm \scriptscriptstyle S}^{\rm \scriptscriptstyle D  }[k]\\
& \quad P_{\rm \scriptscriptstyle e}^{\rm \scriptscriptstyle Purch  } [k] \geq 0\\
& \quad \text{for all } k=0,...,N_{\rm\scriptscriptstyle H}
\end{align*}
where $\xi[k] = [u_{\rm\scriptscriptstyle H}[k],\Delta_{\rm\scriptscriptstyle H}[k], z_{\rm\scriptscriptstyle H}[k] ]'$  and $N_{\rm\scriptscriptstyle H}=24\times 4$ is the number of sampling times included in the $24$-hours horizon and $q_{\rm \scriptscriptstyle S}^{\rm \scriptscriptstyle D }$ is the steam demand. The high level optimization generates the profiles of binary and continuous inputs that minimize the economic cost function of the generation unit for the given utility demand.
\subsection{The low level control}
\label{sec:control-sub:LL}
The low level control, presented in this section, aims to regulate the system to the operating points prescribed by the UC layer.
In this section, we focus on the dynamic control of the boiler unit, as it must satisfy input and output constraints to operate safely, as described in Section \ref{chap:model}. The CHP unit has to just operate at a fixed rotational speed $\bar{\omega} = 50$ Hz and physical constraints are not imposed, thus a standard PI controller is in place. On the other hand, the low level control of the FTB is critical. In particular, here we consider the control problem for the boiler in full operation modes $\rm{ON}_{u_{\rm\scriptscriptstyle ex}=0}$ and $\rm{ON}_{u_{\rm\scriptscriptstyle ex}=1}$, while the proposal of innovative control strategies for the start-up mode will be discussed in Section \ref{chap:ltv_mpc}.\\
In operating condition, the water level $l_{\rm\scriptscriptstyle w}$ and the steam pressure $p_{\rm\scriptscriptstyle s}$ must be maintained in the bounded ranges \eqref{eq:l_const} and \eqref{eq:p_const}. Violating the lower bound in \eqref{eq:l_const} would expose the boiler tubes outside the water, with their subsequent overheating and dangerous consequences while, if the upper bound is violated, water or wet steam may enter the distribution network. Concerning the pressure, the constraints are defined based on the range acceptable by the steam consumers: the interval dimension is application-dependent and the fallout from constraints violation stems from degradation of the steam-user performance to the breakdown of the downstream systems. Note that the pressure constraints can be straightforwardly converted to constraints on the water temperature $T_{\rm\scriptscriptstyle w}$, as the system operates at saturation.\\
The objective of the low-level controller is to enforce these limitations while steering the boiler outputs to the nominal setpoints: a model predictive control  can manage explicitly the constraints on inputs and outputs to easily regulate MIMO systems.\\
The linear discrete-time model of the boiler~\eqref{eq:nl_RK4} is used by the low level controller to predict the system behaviour. For the control of the nominal operation, in production mode, one fixed nominal linearization point has been selected. The state, input, and output vectors of the state-space model are  $x_{\rm\scriptscriptstyle L}= [ T^{\rm\scriptscriptstyle t,g}, T^{\rm\scriptscriptstyle t,CHP},  l_{\rm\scriptscriptstyle w}, T_{\rm\scriptscriptstyle w} ]'$, $u_{\rm\scriptscriptstyle L}=[q_{\rm \scriptscriptstyle f}, q_{\rm \scriptscriptstyle g}^{\rm \scriptscriptstyle B}]'$, and  $y_{\rm\scriptscriptstyle L}= [l_{\rm\scriptscriptstyle w}, p_{\rm\scriptscriptstyle s}^{\rm \scriptscriptstyle B} ]'$. 
At low-level, we consider the steam flow rate $q_{\rm \scriptscriptstyle s}^{\rm \scriptscriptstyle B}$ and the term $u_{\rm\scriptscriptstyle ex} H^{\rm \scriptscriptstyle CHP}_{\rm\scriptscriptstyle ex}$ as measured disturbances: while the latter is not manipulable (since $H^{\rm \scriptscriptstyle CHP}_{\rm\scriptscriptstyle ex}$ is an output of the CHP and the state of $u_{\rm\scriptscriptstyle ex}$ is selected by the high level optimizer), the former is indeed manipulable, although a gas flow request $q_{\rm \scriptscriptstyle s}^{\rm \scriptscriptstyle D}$ is assumed to be given, to be fulfilled at any time instant. However, this request may not be feasible in some cases, e.g., in transient conditions and when the real demand is different from the forecast used for scheduling at the high level. For this reason, it is beneficial for the low-level control to use $q_{\rm \scriptscriptstyle s}^{\rm \scriptscriptstyle B}$ as a further degree of freedom, to be set equal to the demanded $q_{\rm \scriptscriptstyle s}^{\rm \scriptscriptstyle D}$ if possible.\\
The linearization point is characterized by the values $x_{\rm\scriptscriptstyle L}^{\rm\scriptscriptstyle ss}$, $u_{\rm\scriptscriptstyle L}^{\rm\scriptscriptstyle ss}$, $y_{\rm\scriptscriptstyle L}^{\rm\scriptscriptstyle ss}$, $q_{\rm \scriptscriptstyle s}^{\rm\scriptscriptstyle ss}$, and $H^{\rm \scriptscriptstyle CHP,ss}_{\rm\scriptscriptstyle ex}$, while the vectors of the corresponding linearized system
\begin{align}
\delta x_{\rm\scriptscriptstyle L}[h+1] & = A_{\rm\scriptscriptstyle L} \delta x_{\rm\scriptscriptstyle L} [h]+ B_{\rm\scriptscriptstyle L,u} \delta u_{\rm\scriptscriptstyle L} [h]+ B_{\rm\scriptscriptstyle L,q} \delta q_{\rm\scriptscriptstyle s}[h] + B_{\rm\scriptscriptstyle L,ex} \delta v_{\rm\scriptscriptstyle ex}[h] \label{eq:sys_lin1}\\
\delta y_{\rm\scriptscriptstyle L}[h] & = C_{\rm\scriptscriptstyle L} \delta x_{\rm\scriptscriptstyle L}[h] + D_{\rm\scriptscriptstyle L,u} \delta u_{\rm\scriptscriptstyle L}[h] + D_{\rm\scriptscriptstyle L,q} \delta q_{\rm\scriptscriptstyle s}[h]+D_{\rm\scriptscriptstyle L,ex} \delta v_{\rm\scriptscriptstyle ex}[h] \label{eq:sys_lin2}
\end{align}
are defined as
$\delta x_{\rm\scriptscriptstyle L}=x_{\rm\scriptscriptstyle L}-x_{\rm\scriptscriptstyle L}^{\rm\scriptscriptstyle ss}$, $\delta u_{\rm\scriptscriptstyle L}=u_{\rm\scriptscriptstyle L}-u_{\rm\scriptscriptstyle L}^{\rm\scriptscriptstyle ss}$, $\delta y_{\rm\scriptscriptstyle L}=y_{\rm\scriptscriptstyle L}-y_{\rm\scriptscriptstyle L}^{\rm\scriptscriptstyle ss}$, $\delta q_{\rm\scriptscriptstyle s}=q_{\rm\scriptscriptstyle s}^{\rm \scriptscriptstyle B}-q_{\rm\scriptscriptstyle s}^{\rm\scriptscriptstyle ss}$, $\delta v_{\rm\scriptscriptstyle ex}=\alpha_{\rm\scriptscriptstyle ex} H^{\rm \scriptscriptstyle CHP}_{\rm\scriptscriptstyle ex}-H^{\rm \scriptscriptstyle CHP,ss}_{\rm\scriptscriptstyle ex}$.\\
The main objective of the controller is, if possible, to fulfill the steam demand, $\delta q_{\rm\scriptscriptstyle s}^{D}[h]=q_{\rm\scriptscriptstyle s}^{D}[h]-q_{\rm\scriptscriptstyle s}^{\rm\scriptscriptstyle ss}$, while regulating the output $\delta y_{\rm\scriptscriptstyle L}$ to the nominal operating condition $(0,0)$. The input setpoint is defined as $\delta u_{\rm\scriptscriptstyle L}^{\rm\scriptscriptstyle opt}= u_{\rm\scriptscriptstyle L}^{\rm\scriptscriptstyle ref}-u_{\rm\scriptscriptstyle L}^{\rm\scriptscriptstyle ss}$, where $u_{\rm\scriptscriptstyle L}^{\rm\scriptscriptstyle ref}$ is defined based on $q_{\rm\scriptscriptstyle s}^{D}[h]$ and $H^{\rm \scriptscriptstyle CHP}_{\rm\scriptscriptstyle ex}[h]$ using suitable steady-state maps, e.g. \eqref{eq:Boiler_Stat}. Note that, if the steam demand $q_{\rm\scriptscriptstyle s}^{D}[h]$ and the real CHP exhaust flow rate profile $H^{\rm \scriptscriptstyle CHP}_{\rm\scriptscriptstyle ex}[h]$ are equal to the values assumed at the high scheduling level, the required input setpoint $\delta u_{\rm\scriptscriptstyle L}$ is the nominally optimal one derived by the high level optimizer.\\
In addition, the controller must guarantee the operational constraints to be satisfied.\\
The model predictive controller defines the input by solving the following quadratic programming problem, at each time step~$h$:
\begin{align}
\min_{\delta u_{\rm\scriptscriptstyle L}[h],\dots,\delta u_{\rm\scriptscriptstyle L}[h+p],\delta q_{\rm\scriptscriptstyle s}[h],\dots,\delta q_{\rm\scriptscriptstyle s}[h+p]} & \quad J_{\rm\scriptscriptstyle L}
\end{align}
\begin{align*}
\begin{array}{lrl}
\text{s.t.} & &\quad	{\rm model} \, \eqref{eq:sys_lin1}-\eqref{eq:sys_lin2} \\
& &\quad	\delta x_{\rm\scriptscriptstyle L}[h]= \delta \hat{x}_{\rm\scriptscriptstyle L}[h|h] \\
 & 	\delta q_{\scriptscriptstyle \rm f}\in &\quad [0 - q^{\rm\scriptscriptstyle ss}_{\scriptscriptstyle \rm f},q_{\scriptscriptstyle \rm f,max} - q^{\rm\scriptscriptstyle ss}_{\scriptscriptstyle \rm f}] \\
&	\delta q^{\scriptscriptstyle \rm B}_{\rm\scriptscriptstyle g}\in &\quad[q^{\scriptscriptstyle \rm B}_{\scriptscriptstyle \rm min} - q^{\scriptscriptstyle \rm B,ss}_{\scriptscriptstyle \rm g}, q^{\scriptscriptstyle \rm B}_{\scriptscriptstyle \rm max} - q^{\scriptscriptstyle \rm B,ss}_{\scriptscriptstyle \rm g}]\\
&	\delta T_{\scriptscriptstyle \rm w}\in &\quad[T_{\scriptscriptstyle \rm min} - T_{\scriptscriptstyle \rm w}^{\scriptscriptstyle \rm ss}, T_{\scriptscriptstyle \rm max} - T_{\scriptscriptstyle \rm w}^{\scriptscriptstyle \rm ss}] \\
&	\delta l_{\scriptscriptstyle \rm w}\in&\quad [l_{\scriptscriptstyle \rm min} - l_{\scriptscriptstyle \rm w}^{\scriptscriptstyle \rm ss}, l_{\scriptscriptstyle \rm max} - l_{\scriptscriptstyle \rm w}^{\scriptscriptstyle \rm ss}]
\end{array}
\end{align*}
where
$J_{\rm\scriptscriptstyle L}=\sum_{j=h}^{h+p} \|\delta y_{\rm\scriptscriptstyle L}[j]\|_{W_{Q}}^2 + \|\delta u_{\rm\scriptscriptstyle L}[j] - \delta u_{\rm\scriptscriptstyle L}^{\rm\scriptscriptstyle opt}[j]\|_{W_{R}}^2+\|\delta q_{\rm\scriptscriptstyle s}[j] - \delta q_{\rm\scriptscriptstyle s}^{\rm\scriptscriptstyle D}[j]\|_{W_{Q_{\rm\scriptscriptstyle s}}}^2$
and where $\delta \hat{x}_{\rm\scriptscriptstyle L}[h|h]$ is the current state estimate obtained by a proper observer, see Section~\ref{sub:obs}. A slack variable can be introduced to further enforce the feasibility of the optimization problem at all time instants. The matrices $W_{Q}$, $W_{R}$, $W_{Q_{\rm\scriptscriptstyle s}}$ are properly-defined positive-definite matrices; note that, in order to give priority to the fulfilment of the steam demand, we can set $W_{Q_{\rm\scriptscriptstyle s}}>>\lambda_{\rm\scriptscriptstyle max}(W_{Q}),\lambda_{\rm\scriptscriptstyle  max}(W_{R})$.\\
At time instant $h$, the optimal values $\delta u_{\rm\scriptscriptstyle L}(h|h)$ and $\delta q_{\rm\scriptscriptstyle s}(h|h)$ are obtained: therefore, the input $u_{\rm\scriptscriptstyle L}[h]= u_{\rm\scriptscriptstyle L}^{\scriptscriptstyle \rm ss}+\delta u_{\rm\scriptscriptstyle L}(h|h)$ is applied to the real system and the steam flow rate $q_{\rm\scriptscriptstyle s}^{\rm \scriptscriptstyle B}[h]= q_{\rm\scriptscriptstyle s}^{\scriptscriptstyle \rm ss}+\delta q_{\rm\scriptscriptstyle s}(h|h)$
is actually provided.\\
It is important to note that, if the steam demand $q_{\rm\scriptscriptstyle s}^{\rm \scriptscriptstyle D}[h]$ and/or the real CHP exhaust flow rate profile $H^{\rm \scriptscriptstyle CHP}_{\rm\scriptscriptstyle ex}[h]$ differ from the values considered at the high scheduling level, the steady-state points reached by the controlled system will possibly differ from the ones selected by the high level optimizer. This, however, is not critical since, from the practical point of view, the major goal of the controller is to guarantee that the operational constraints are enforced at all time instants.\\
	The steady-state quantities, e.g. $H^{\rm \scriptscriptstyle CHP,ss}_{\rm\scriptscriptstyle ex}$ and $q_{\rm\scriptscriptstyle s}^{\rm\scriptscriptstyle ss}$, are defined based on the high-level optimization and depend on the forecast of the demand considered at high level.
	When the optimization is executed, it typically considers a day-ahead prediction of the demand over a 24-hour horizon, discretized by time step $\tau_{\rm H}$. This is in general a piece-wise constant approximation of the  demand forecast. During operation, more accurate predictions are available, as well as a measure of the current/past demand mismatch with respect to the profile considered in the optimization.\\
	If, in any case, the actual values of $q_{\rm\scriptscriptstyle s}^{\rm \scriptscriptstyle D}[h]$ and/or $H^{\rm \scriptscriptstyle CHP}_{\rm\scriptscriptstyle ex}[h]$ differ dramatically  from the ones considered for the high-level optimization, the top-layer optimizer may be required to be run again (with new and more reliable forecast of steam and electric power demands) in order to recompute, for the following time steps, more realistic optimal set-points.
	We define thresholds $\varepsilon_H$ (or $\varepsilon_q$) that determine the maximum acceptable mismatch, $|H^{\rm \scriptscriptstyle CHP}_{\rm\scriptscriptstyle ex}[h] -H^{\rm \scriptscriptstyle CHP,ss}_{\rm\scriptscriptstyle ex}|< \varepsilon_H$ (or $|q_{\rm\scriptscriptstyle s}^{\rm \scriptscriptstyle D}[h] - q_{\rm\scriptscriptstyle s}^{\rm\scriptscriptstyle ss}| < \varepsilon_q$), otherwise   the computation of a new optimization at high level is triggered.
	In this scenario, it can be reasonable to use updated forecast of the future demand.	

\subsubsection{Stand-by mode} \label{sec:stand_by}
As discussed, when the consumer demand  is lower than the minimum level, i.e., when $q_{\rm \scriptscriptstyle s}^ {\rm \scriptscriptstyle D} < \bar{q}_{\rm \scriptscriptstyle s}^ {\scriptscriptstyle \min}$, the FTB can be either shut down or put in the stand-by mode.
\begin{figure}[htb]
	\centering
	\includegraphics[width=1\linewidth]{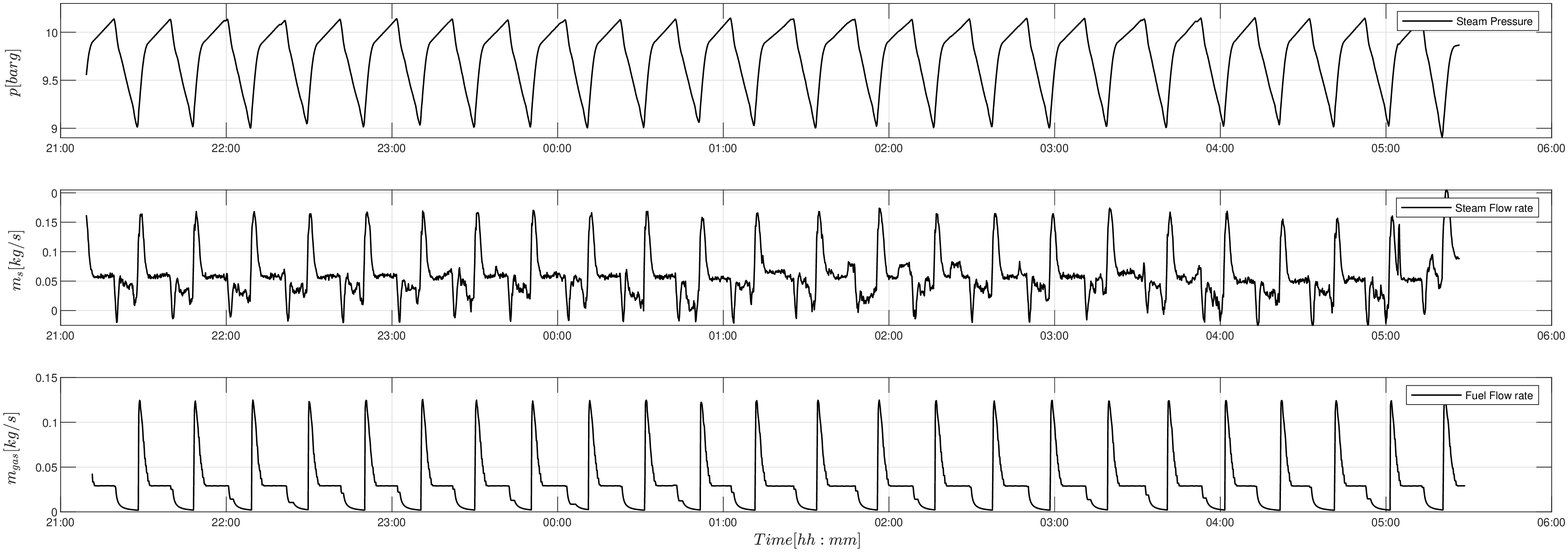}
	\caption{Limit cycle operation in stand-by mode}
	\label{fig:Cyclic_dyn}
\end{figure}%
The stand-by mode, activated by the boolean input $u_{\rm\scriptscriptstyle SB}^{\rm\scriptscriptstyle B}$, 
is a low-power setting operation, characterized by a limit cycle, where the steam pressure is controlled to lie in a broad range and only a minimum amount of steam is flowing to keep the distribution network warm. In this mode, the FTB can be immediately restarted to follow the required steam demand, as soon as the latter is restored above the minimum threshold.\\
Figure \ref{fig:Cyclic_dyn} shows available experimental data, where the limit cycle is characterized by a low frequency fluctuation of the pressure - with a period of about 15 minutes - confined in the range $[\bar{p}_{\scriptscriptstyle \rm min}, \bar{p}_{\scriptscriptstyle \rm max}]$.\\ Pressure thresholds 
trigger a simple three-mode controller used to regulate the system maintaining a low natural gas consumption.\\
\begin{itemize}
	\item \textit{mode 1:} as the pressure reaches its low limit, $\bar{p}_{\scriptscriptstyle \rm min}$, the burner abruptly provide a high heat impulse until nominal pressure level, $\bar{p}_{\scriptscriptstyle \rm sp}$, is achieved;
	\item \textit{mode 2:} the controller keeps the fuel flow rate steady to a low value, reducing the pressure rise until the maximum value, $\bar{p}_{\scriptscriptstyle \rm max}$, is reached;
	\item \textit{mode 3:} the burner is set at the lowest firing mode, which causes a degradation of the system pressure.
\end{itemize}
\begin{figure}[thpb]
	\centering
	\includegraphics[width=.6\linewidth]{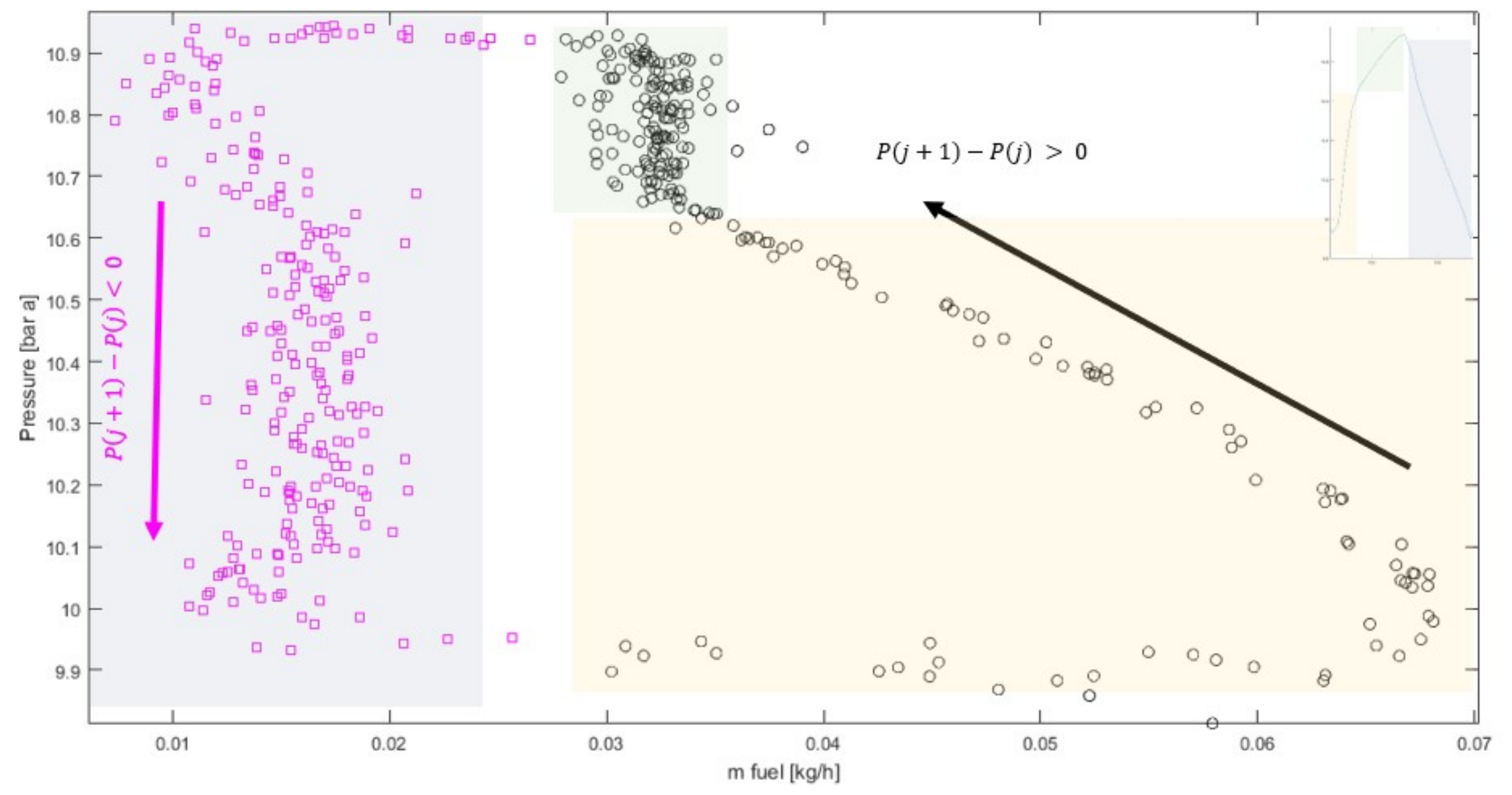}
	\caption{Three-mode controller used in stand-by mode: scatter plot of the pairs pressure-fuel flow rate for rising (black circles) and decreasing pressure (magenta squares). Three regions, with different background colors, show the controller modes: mode 1 in light yellow, mode 2 in green, mode 3 in grey. In the up-right corner the pressure cycle}
	\label{fig:Cyclic_Controller}
\end{figure}%
The controller strategy is identified from experimental data analyzing the relation between absolute pressure and fuel flow rate, as show in  Figure  \ref{fig:Cyclic_Controller}. In the up-right corner the pressure cycle is reported for clarity: the three regions are highlighted by different shaded background patches, while pressure-fuel flow rate pairs are shown in black circles, for rising pressure and in magenta squares, for decreasing pressures.\\
A PI controller is formulated to simulate the fuel profile for the first phase of the cycle:
\begin{align*}
q_{\scriptscriptstyle \rm g}[h] = q_{\scriptscriptstyle \rm g}^{\scriptscriptstyle \rm 0} + K_{\scriptscriptstyle \rm P} (p[h] -\bar{p}_{\scriptscriptstyle \rm PI}) + K_{\scriptscriptstyle \rm I} \sum_{\scriptscriptstyle j=h'}^{\scriptscriptstyle h} (p[j]-\bar{p}_{\scriptscriptstyle \rm PI})
\end{align*}%
where the controller parameters, $K_{\scriptscriptstyle \rm P} = 2.77, K_{\scriptscriptstyle \rm I}=0.014, q_{\scriptscriptstyle \rm g}^{\scriptscriptstyle \rm 0}= 0.066,\bar{p}_{\scriptscriptstyle \rm PI}$ = 9.7, are identified from data and $h'$ is the switching time from phase 3 to phase 1, when the integral error is reset to zero.
In the other phases, a steady value of gas flow rate is defined.
\subsection{State observer}\label{sub:obs}
In this paper, we consider a configuration of the system in which the state vector is not fully accessible, in contrast with the previous work \cite{Spinelli2019}: in this preliminary work, the optimal startup procedure is  investigated for a simpler boiler model with the main assumption of full state available.\\
	 In this work, this assumption is removed in order to fit in a more realistic scenario. In fact, here both the regulating MPC, for the production modes, and the LPV-MPC, defined for the optimization of the start-up procedure, rely on the actual sensor existing on the system in their real locations. \\
	As described in Section \ref{sub:boiler_low}, the fire tube boiler is endowed with a meter for the measure of the water level and an internal pressure sensor.
An observer is set up to reconstruct the state vector based on the available measurements. A discrete-time Extended Kalman Filter has been designed to this aim, as discussed e.g. in \cite{Simon2006}.
We assume that a process and measurement noises, denoted respectively $w_{\scriptscriptstyle \rm pr}[h]$ and $w_{\scriptscriptstyle \rm m}[h]$, act on the discrete time system derived in \ref{sub:boiler_low}:
\begin{align*}
\label{eq:nl_RK4_2}
x[h+1] = & f_{RK4}(x[h],u[h]) + w_{\scriptscriptstyle \rm pr}[h]\\
y[h] =  & {g}(x[h],u[h]) + w_{\scriptscriptstyle \rm m}[h]
\end{align*}
where $w_{\scriptscriptstyle \rm pr}[h]$ and $w_{\scriptscriptstyle \rm m}[h]$ are zero mean Gaussian noises with covariances $X_{Q}$ and $X_{R}$, compatible with the quantified measurement errors. Their values are reported in Section \ref{sec:opt_res}. Here the dependency of the model upon the parameter vector $\theta$ is discarded for simplicity of notation.

\section{Optimal start-up procedure}
\label{chap:ltv_mpc}
The control of the start-up mode differs dramatically from the dynamic control of the production mode, regarding the objective and the constraints involved, but most of all due to the increased importance of the system non-linearity, which precludes the direct extension of the low-level control to the management of the start-up procedure.\\
We propose a nonlinear MPC approach for the boiler start-up. The proposed approach formulates the MPC problem including the additional optimization variable ${x}^{\rm \scriptscriptstyle ss}_{ \scriptscriptstyle h}$, which represents a ''temporary'' target steady state, reachable in $N_{\rm \scriptscriptstyle P}$ discrete-time steps (i.e., at the end of the optimization horizon), as the closest one with respect to the desired final state. A further constraint is introduced to force the new decision variable ${x}^{\rm \scriptscriptstyle ss}_{ \scriptscriptstyle h}$ to be an intermediate equilibrium point. This can guarantee practically the recursive feasibility. A formal proof of this feature is left for future research work.\\
This intermediate steady state is used as a terminal constraint and as a reference value in the objective function. Thanks to this, the optimization horizon can be greatly decreased, therefore reducing the numerical complexity, and, as a byproduct, to avoid the (often complex) computation of a terminal positively invariant set and of a suitable terminal cost and to guarantee recursive feasibility of the MPC optimization problem.
Moreover, this approach does not require to know the actual minimal duration of the start-up phase.\\
This method, by considering only a small percentage of the entire start-up duration in its prediction horizon, compromises the global optimality of the solution. However, this suboptimality is mitigated by considering in the cost function a penalization of the distance of the temporary steady-state with respect to the final nominal target, according to the dynamic programming paradigm. The level of this performance degradation will be evaluated numerically on the case study.\smallskip\\
To achieve an efficient implementation of the NMPC, we propose here a parameter-varying linearisation of the nonlinear system, computed along the state/input trajectory obtained at the precedent optimization cycle. This reduces the nonlinear optimization to a Constrained Quadratic Program (QP).
This approach is referred here to as Linear Parameter-Varying Model Predictive Control (LPV-MPC) and it is similar to the one proposed in \cite{Falcone08}. The main difference with \cite{Falcone08} consists of how the trajectory around which the system is linearised is computed. A strong connection with the Real-Time iteration scheme proposed in \cite{Diehl2005} also exists.
As a more general formulation of model predictive control, LPV-MPC can be directly used both in start-up optimal control and in the regulation of productive modes. However the boiler stabilization in nominal operating condition is achievable in the most efficient way by switching to a simpler linear model predictive control, where the prediction model is obtained by linearization at the nominal operating point.
\begin{figure}[tb]
	\centering	
	\includegraphics[width=0.55\linewidth]{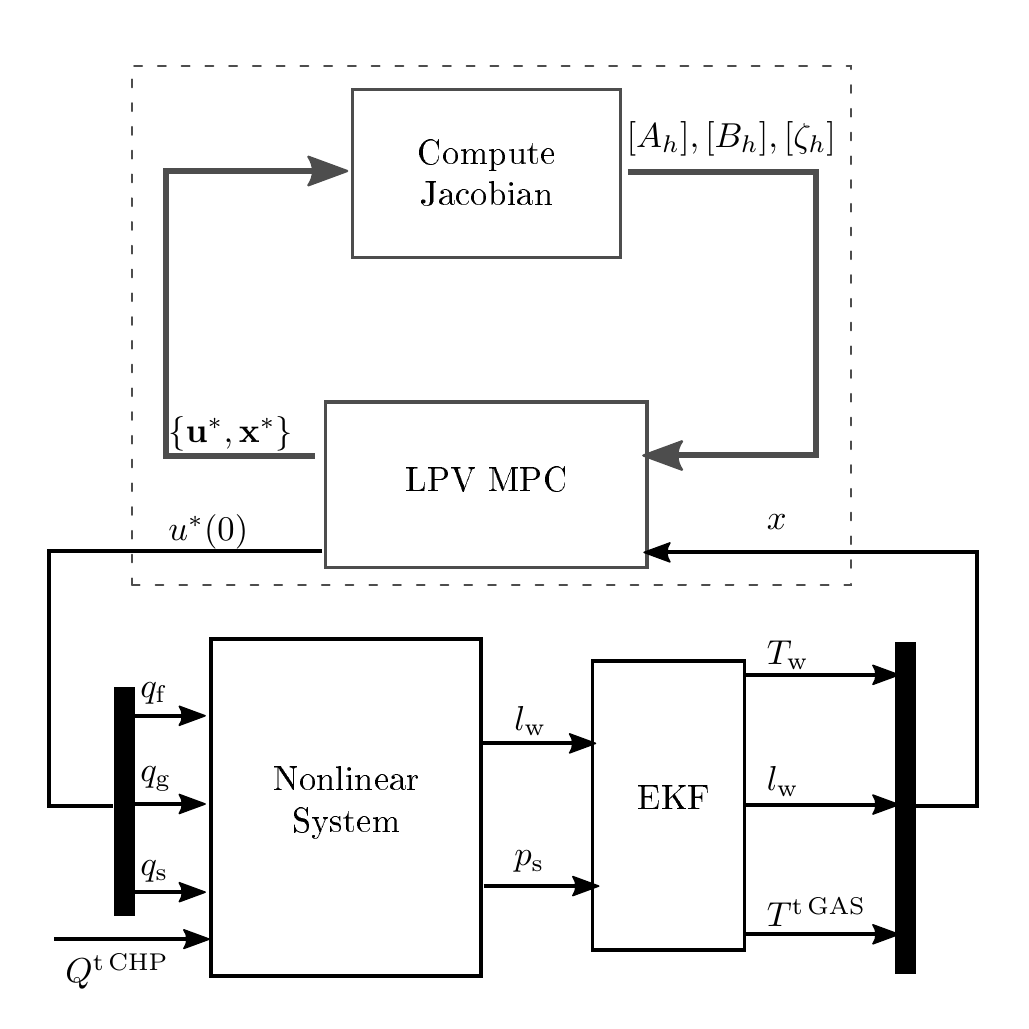}
	\caption{Linear Parameter-Varying MPC scheme of the FTB system}
	\label{fig:LPV_MPC_FTB}
\end{figure}%
The system~\eqref{eq:nl_model_tot} is integrated with the sampling time $\tau_{\rm\scriptscriptstyle S}=6$ s for better numerical accuracy. For control purposes, at each time step (say $h$), the discretized model is linearized around the last available optimal state/input trajectories, e.g., $x^{\scriptscriptstyle}[i|h-1]$ and $u^{\scriptscriptstyle}[i|h-1]$, respectively, obtained at step $h-1$ is, for $i=h,\dots,N_p-1$
\begin{equation}
\label{eq:l_RK4}
x[i+1] = A_{i|h-1} x[i] + B_{i|h-1} u[i] + \zeta_{i|h-1}
\end{equation}
where
\begin{align*}
A_{i|h-1} = &\frac{\partial f_{RK4}}{\partial x}\bigg\rvert_{x^{\scriptscriptstyle}[i|h-1],u^{\scriptscriptstyle}[i|h-1]} \quad B_{i|h-1} = \frac{\partial f_{RK4}}{\partial u}\bigg\rvert_{x^{\scriptscriptstyle}[i|h-1],u^{\scriptscriptstyle}[i|h-1]}\\
\zeta_{i|h-1} = & f_{RK4}(x^{\scriptscriptstyle}[i|h-1],u^{\scriptscriptstyle}[i|h-1]) - A_{i|h-1} x^{\scriptscriptstyle}[i|h-1] - B_{i|h-1} u^{\scriptscriptstyle}[i|h-1]
\end{align*}
Given the final state target $x_{T}$, the optimization program to be solved at time step $h$ is the following
  \begin{equation}
\label{eq:LPV-MPC}
\begin{split}
\min_{u[h],\dots,u[h+N_p-1], {x}^{\rm \scriptscriptstyle ss}_{ \scriptscriptstyle h}} & \sum_{i=h}^{h+N_p} \norm{x[i] - {x}^{\rm \scriptscriptstyle ss}_{ \scriptscriptstyle h}}^2_{W_{Q}} + \norm{u[i]-u[i-1]}^2_{W_{R}} + \norm{{x}^{\rm \scriptscriptstyle ss}_{ \scriptscriptstyle h} - x_{T}}^2_{W_{Q_T}}\\
s.t. & \quad \text{dynamics }\eqref{eq:l_RK4}\\
& \quad x[h] = \hat{x}[h|h] \\
&\quad {c}(x[i],u[i]) \leq 0 \,\,\text{for all }i=h,\dots,h+N_p-1\\
& \quad x_j[i+1]-x_j[i] \leq  \sigma^{Max}\tau_{\rm\scriptscriptstyle S} \quad \text{for all } j \in \mathcal{J}_{\sigma}\text{ and for all }i=h,\dots,h+N_p-1\\
& \quad x[h+N_p] = {x}^{\rm \scriptscriptstyle ss}_{ \scriptscriptstyle h}\\
&\quad {c}({x}^{\rm \scriptscriptstyle ss}_{ \scriptscriptstyle h},{u}^{\rm \scriptscriptstyle ss}_{ \scriptscriptstyle h}) \leq 0
\end{split}
\end{equation}
In the latter problem, $\hat{x}[h|h]$ is the estimated state obtained using the state estimator described in Section~\ref{sub:obs} (designed at sampling time $\tau_{\rm\scriptscriptstyle S}$). Also, function $c(\cdot,\cdot)$ includes the constraints \eqref{eq:boiler_CONSTRAINT}. These constraints are also applied the terminal steady state condition ${x}^{\rm \scriptscriptstyle ss}_{ \scriptscriptstyle h}$; ${u}^{\rm \scriptscriptstyle ss}_{ \scriptscriptstyle h}$ is the corresponding steady-state input. Finally, a hard constraint is enforced to limit the thermal stress consistently with~\eqref{eq:th_bound}: it is applied to each component $j$ of the state vector that are correlated to thermal stress, included for simplicity in the set $\mathcal{J}_{\sigma}$.\\
The choice of weighting matrices is not a trivial task. The weight matrices used in the stage cost are usually diagonal matrices, here defined by considering the inverse of the squares of the nominal values for each state and input channel, thus normalizing the absolute contribution of each channel. This can simplify the tuning of the diagonal entries of $W_Q$ and $W_R$, that can reflect the preferences of the controller designer.\\
For what concern the last term of the cost function, measuring the distance from the actual target, the idea for tuning the weighting matrix is to emulate the approaches presented for designing the terminal cost of MPC.\\ One of the classic approaches is based on the computation of an infinity horizon LQR regulator, i.e., by considering the solution $P$ of a discrete algebraic Riccati equation associated to the system linearized at the end of the horizon $A_{N_{\rm{P}-1}|h-1}$, $B_{N_{\rm{P}-1}|h-1}$, with the given $W_Q$ and $W_R$.\\ Thus, $W_{Q_T}$ can be selected to be $W_{Q_T}=\alpha P$, where $\alpha$ is a designer tuning parameter.%

\section{Simulation results} \label{chap:sim_res}
In the following sections we validate the proposed hierarchical optimization approach and the start-up optimal control through simulation, using the nonlinear models of the GU presented in Section \ref{chap:model}. The water and steam thermophysical properties, based on IAPWS-IF97, are computed using \cite{SteamTool}.\\
The MLD model \eqref{eq:MLD} is implemented using the hybrid system description language (HYSDEL) \cite{HYS2004}, where transition logics and switched affine systems are formulated directly using mixed-integer inequalities \cite{Mitra94,Williams2013}.\\
The optimal control is implemented in Matlab using YALMIP, with Gurobi 8.1.0 as a solver for the MILP and QP problems. Simulations are performed on PC with a quad-core processor 
and 16GB RAM.
 The closed-loop simulations are performed using the nonlinear models presented in Section \ref{chap:model}, for the boiler and the CHP. As concerning the latter, the nonlinear model implemented is reported in the original papers \cite{Guzzella2004} and \cite{Videla2007}.
	The continuous-time ODEs are integrated using variable-step variable-order methods for stiff ODEs (as implemented in Matlab ode15s routine) considering a interval of integration defined by the controller sampling time.
	The models used for control are linearized model obtained from the nonlinear models by a prior discretization, using a Runge-Kutta method of the 4th order.

\subsection{Hierarchical optimization}
\begin{figure}[tb]
	\centering
	\includegraphics[width=1\linewidth]{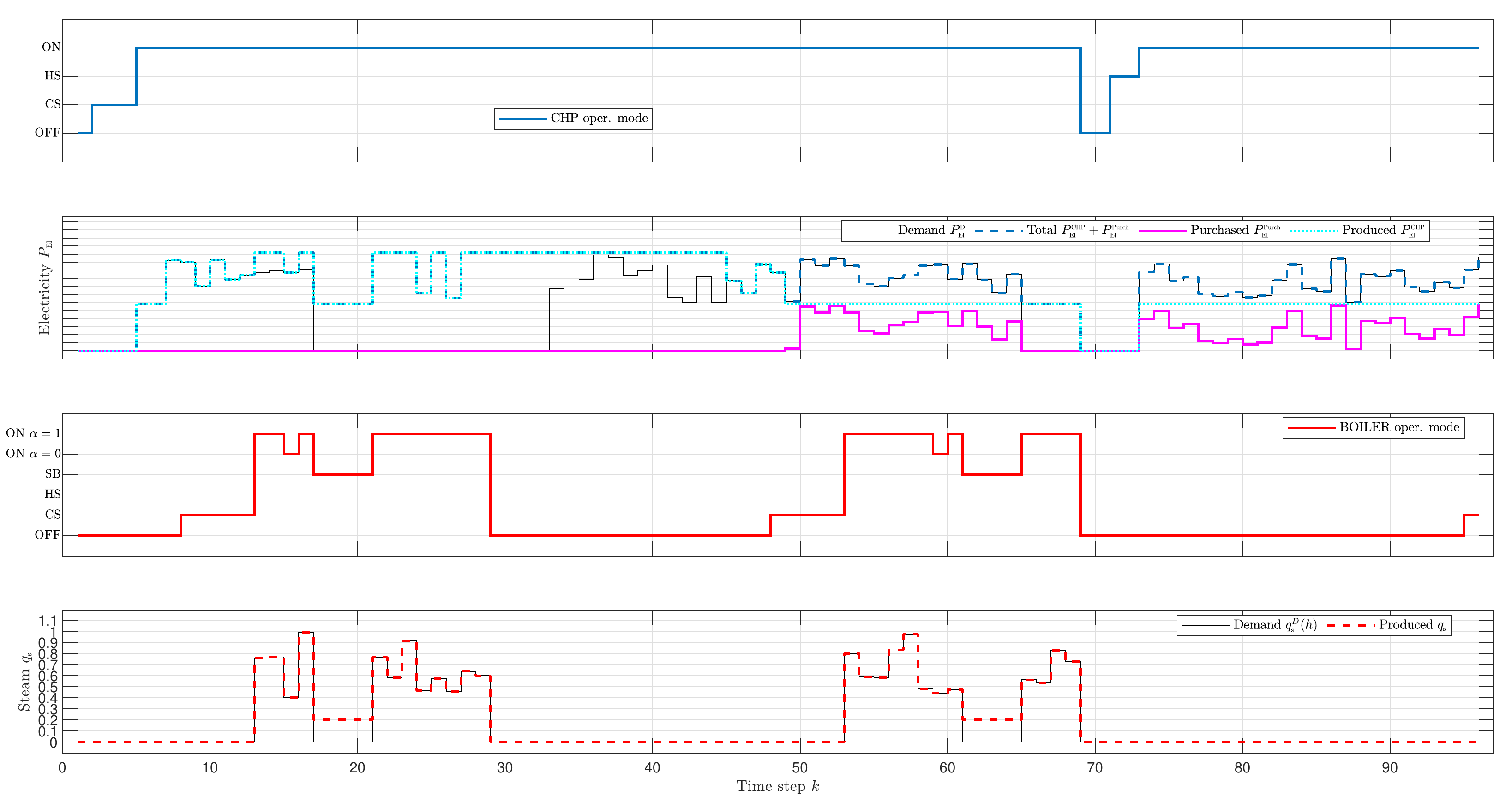}
	\caption{Simulation results of high level optimization: at $\tau = 48$ the gas price is switched to high. i) CHP operating modes. ii) Electricity demand and supply. iii) FTB operating modes. iv) Steam demand and production.}
	\label{fig:Sim}
\end{figure}

First, we present the results on the hierarchical control scheme for a network composed by a FTB and a CHP for electricity production.
The steam generator is a three-pass fire tube boiler working at a nominal pressure of $10$~bar and with maximum steam flow-rate of $\bar{q}_s = 1200$~{kg}/{h}, while the CHP is a $12$~valve natural gas ICE producing up to $\bar{P}_{\rm \scriptscriptstyle e}^{\rm \scriptscriptstyle CHP} =1200$~kW.
The FTB and the ICE are governed to supply the demanded steam and electrical power, operating within the permitted ranges summarized in Table~\ref{table_example}.\\
\begin{table}[tb]
	\caption{Lower and Upper Bounds on the FTB and CHP Variables}
	\label{table_example}
	\begin{center}
		\begin{tabular}{|l||l||l|}
			\hline
			Variable & Minimum & Maximum\\
			\hline
			$P_{\rm \scriptscriptstyle e}^{\rm \scriptscriptstyle CHP}$ &  $50\%$ & $100\%$\\
			\hline
			$q^{\scriptscriptstyle \rm B}_{\rm\scriptscriptstyle g}$ & $12.5\%$ & $100\%$\\
			\hline
			$q^{\scriptscriptstyle \rm B}_{\rm\scriptscriptstyle f}$ & $0$ & $0.35$~{kg}/{s}\\
			\hline
			$q^{\scriptscriptstyle \rm B}_{\rm\scriptscriptstyle s}$ & $0$ & $0.35$~{kg}/{s}\\
			\hline
			$p_{\rm\scriptscriptstyle s}$ & $ 9.5$~bar & $10.5$~bar\\
			\hline
			$l_{\rm\scriptscriptstyle w}$ & $\bar{l}_{\rm\scriptscriptstyle w} - 0.5\%$ & $\bar{l}_{\rm\scriptscriptstyle w} + 0.5\%$\\
			\hline
			$\sigma_j$ &  $0$ & $6$\\
			\hline
		\end{tabular}
	\end{center}
\end{table}%
\begin{figure}[tbh]
	\centering	
	\includegraphics[width=.55\linewidth]{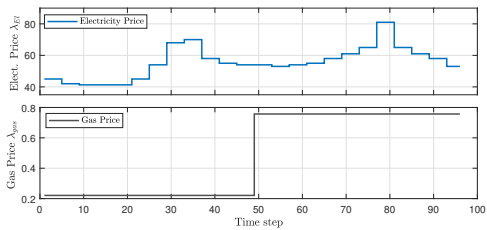}
	\caption{Electricity price fluctuation [\euro$/kWh$] (top). Gas price [\euro$/m^3$] (bottom). At $k=48$ gas price is increased to show by simulation the effect of relative electricity/gas price on the GU behaviour. }
	\label{fig:prices}
\end{figure}
The electricity price is supposed to vary on hourly basis, see Figure \ref{fig:prices}, and its values are taken from the historical database of day-ahead market of \textit{Gestore Mercati Energetici} (GME), the Italian company managing the energy market.\\ For the purpose of highlighting the effect of a different utility price ratio and to show this in a unique graph, the price of natural gas is artificially increased at time $k =48$. A random realization of the fluctuating profiles of electricity and steam demands of the consumers is generated for the simulation.
For the higher layer, mode transition times are reported in Table \ref{table:tau}.\\
\begin{table}[tb]
	\caption{Mode transition times of  FTB and CHP DHA models}
	\label{table:tau}
	\begin{center}
		\begin{tabular}{|l||l||l|}
			\hline
			$\tau$ &Boiler & CHP\\
			\hline
			$\tau_{\scriptscriptstyle \rm{OFF}} $ &$1.0$ & $1.0$\\
			\hline
			$\tau_{\scriptscriptstyle \rm{OFF}\rightarrow \rm{CS}} $ & $5.0$ & $3.0$\\
			\hline
			$\tau_{\scriptscriptstyle \rm{CS}\rightarrow \rm{ON}} $& $4.0$ & $2.0$\\
			\hline
			$\tau_{\scriptscriptstyle \rm{HS}\rightarrow \rm{ON}}$ & $1.0$& $1.0$\\
			\hline
			$\tau_{\scriptscriptstyle \rm{ON}} $ &$3.0$& $3.0$\\
			\hline
			$\tau_{\scriptscriptstyle \rm{SB}\rightarrow \rm{ON}} $ &$0$& -\\
			\hline
		\end{tabular}
	\end{center}
\end{table}
Figure \ref{fig:Sim} shows the results of the high level optimization. The first and the third panels from the top describe the evolution of the operative modes of CHP and FTB, respectively, while the other two panels show the demanded and produced electricity and steam, respectively. The electricity supply is split between the contribution of CHP generation (in dotted blue line) and electricity purchased from the grid (solid magenta). The gas/electricity price ratio affects the optimal strategy of the GU: for $k < 48$, where grid electricity is relatively more expensive, the CHP fulfils the electrical demand alone and, whenever it is convenient, it sells the electricity surplus. In the opposite case, for $k \geq 48$, the CHP is maintained in minimum production mode, just to provide the diverted gas into the FTB for steam production, while the remaining electricity demand is covered purchasing power from the grid.\\
As discussed in Section \ref{sec:stand_by}, when steam demand is zero, the FTB can be either  off or set in stand-by mode. 
The choice between the two strategies depends on the actual duration of the demand interruption: the break-even point, which discriminate the optimal strategy, depends on the downtime and the average consumption in stand-by and in start-up modes.\\
Historical data analysis has shown a reduction of more than $60\%$ of natural gas usage if the stop-and-start strategy is preferred to the stand-by mode, for a no-steam demand of 5 hours. Clearly, a forecast of the future demand is required to optimally decide the best operation mode sequence.
In industrial scenarios, production scheduling can provide useful information to build a proper forecast of the steam usage.
The demand profile in Figure \ref{fig:Sim} is characterized by two steam requests, interspersed with short and long downtime periods, showing different strategies during the short pauses, when FTB is switched to stand-by mode, and the longer periods with zero demand.\\
The reference trajectory computed by the optimization layer is provided to the low-level MPC controller described in Section \ref{sec:control-sub:LL}, which operates when the FTB is in production modes. 
\begin{figure}[t!]
	\centering	
	\includegraphics[width=.95\linewidth]{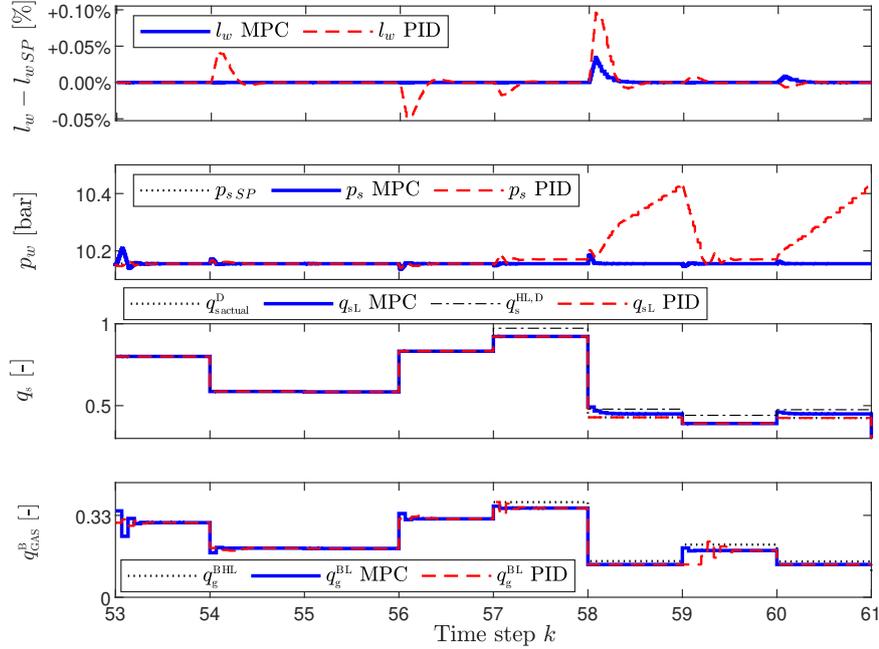}
	\caption{Low-level MPC (solid blue) and PID (dashed red) controllers. 
		i) Water level. ii) Steam pressure. iii) Produced steam flow rate. iv) Boiler gas flow rate.
		\textit{Nominal}: $k =[53, 57)$;  \textit{Scenario $a$}:  $k =[57,61]$, $q_{\rm \scriptscriptstyle s} = q_{\rm \scriptscriptstyle s\,HL} - \delta q_{\rm \scriptscriptstyle s}$.}
	\label{fig:low}
\end{figure}
The linear model is computed around a single nominal operating point of $q_{s\, 0} = 0.5 \bar{q}_s$ and $H_{\rm \scriptscriptstyle ex\, 0}^{\rm \scriptscriptstyle CHP} = 0.5\bar{H}_{\rm \scriptscriptstyle ex}^{\rm \scriptscriptstyle CHP}$.
The MPC regulates the $\delta q_{\rm \scriptscriptstyle GAS}^{\rm \scriptscriptstyle B}$ and $\delta q_{\scriptscriptstyle \rm f}$ to maintain the steam pressure and the water level close to their set points and inside the valid ranges, defined in Table \ref{table_example}.
The advantages of a MIMO optimal controller at low level is highlighted, through a comparison with a classic PID-based control scheme, based on two separate control loops for level and pressure regulation. It shows the robustness of the proposed hierarchical approach, by simulating different scenarios in which the demand used for computing the low-level input set-points differs with respect to the one calculated by the high level optimizer. 
In Figure \ref{fig:low} a simulation displays the behavior of the low level controllers (the MPC in solid blue line and the PID, as a red dashed line), where the nominal steam demand computed at the high layer is considered in the period $k=[53, 57)$, while a disturbed steam demand - lower than the one demanded by the high layer - is considered in the second half of the simulation, where $k=[57, 61)$.
While the MPC regulates in a fast way the system in both nominal and disturbed steam demand, the PI controller - tuned to have comparable performances in the nominal case - can correctly govern the level, but with degradation in pressure regulation.

\subsection{Optimal start-up}\label{sec:opt_res}
A discrete time model is obtained from \eqref{eq:nl_RK4}, with sampling time $\tau_{\rm \scriptscriptstyle S}$.
The sensitivity matrices are calculated by differentiation: symbolic calculations are used to pre-compute off-line the Jacobians $A_{\scriptscriptstyle h}$, $B_{\scriptscriptstyle h}$, based on Automatic Differentiation tools.\\
The start-up optimization is done by solving the QP problem \eqref{eq:LPV-MPC}, with a prediction horizon $N_p=50$ and the following weighting matrices: $W_{Q} = \text{diag}(0.1, 5, 20) \circ w_{x}$, $W_{R} =  \text{diag}(0.01, 0, 0) \circ w_{u}$ and $ W_{Q_T} = \text{diag}(0.1, 5, 30) \circ w_{x}$,
where $\circ$ denotes the Hadamard product, while $w_{x}$, $w_{u}$ contain the inverse of the square of the maximum values of state and input components.\\
The Extended Kalman filter has been implemented considering the following covariance matrices for the noises: 
$X_{Q} = \text{diag}(10^{-2}, 10^{-2}, 10^{-2})$ and $X_{R} = \text{diag}(10^{-3}, 10^{-3})$ and with initial covariance estimate $X_{P \, \scriptscriptstyle 0|0} = [1.5, 10^{-3}, 0.15;10^{-5}, 10^{-3}, 10^{-5};0.15, 10^{-3}, 0.15]$. The state estimate $\hat{x}_{ \scriptscriptstyle 0|0}$ is initialized with the measured level and assigning to both the tube and water temperatures a prescribed initial value.\\

\begin{figure}[!tbh]
	\centering
	\begin{subfigure}[t]{.48\textwidth}
		\centering
		\includegraphics[width=1\linewidth]{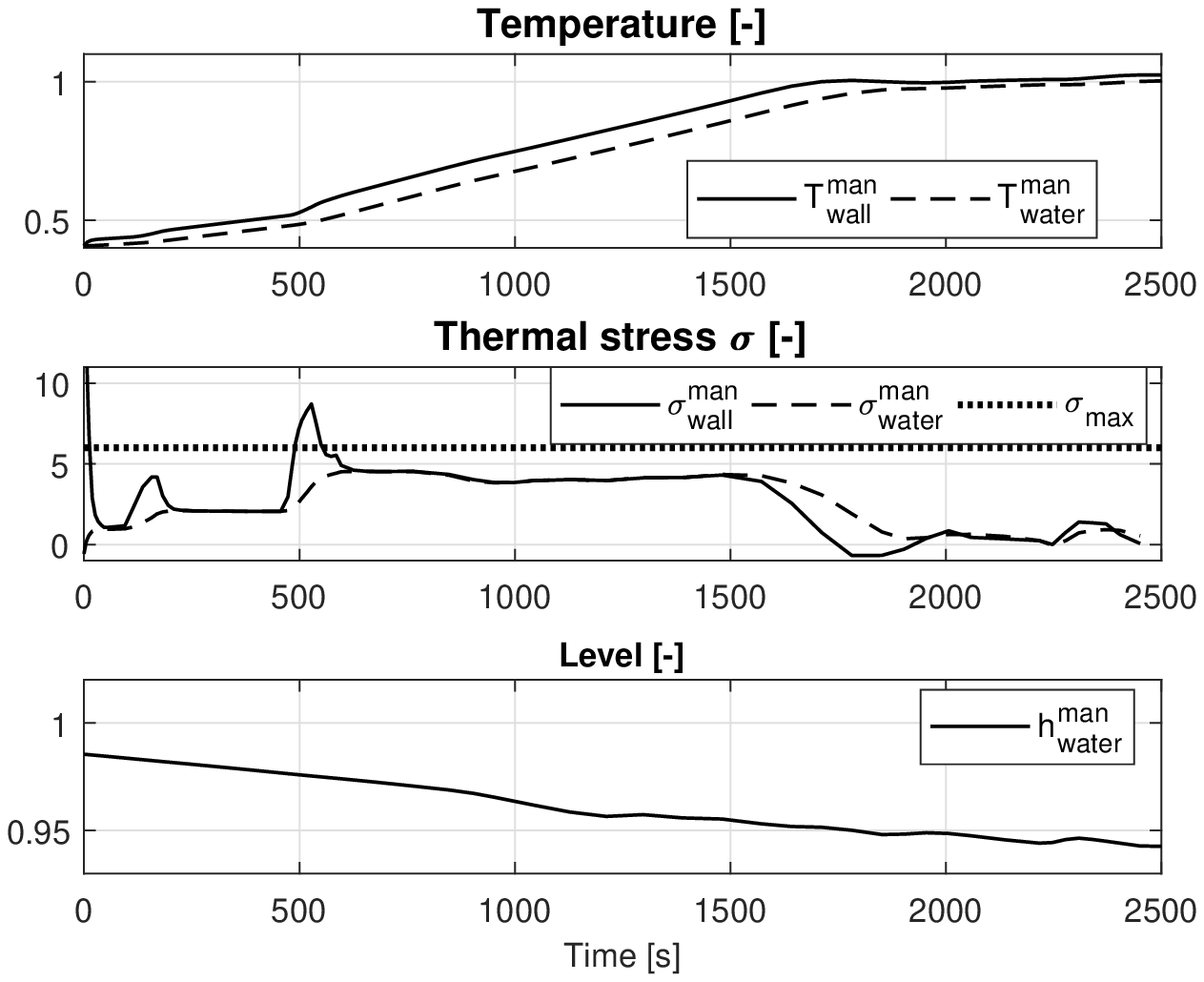}
		\caption{FTB start-up manual procedure - output variables. In the top panel: temperature of fire tubes (solid) and the water temperature (dashed). In the second, the thermal stress on tubes (solid), on the shell (dashed) and its maximum value (dotted). In the third, water level. All graphs are adimensionalized with nominal values, for confidentiality.}
		\label{fig:man_out}
	\end{subfigure}%
	\hfill
	\begin{subfigure}[t]{0.48\textwidth}
		\centering
		\includegraphics[width=1\linewidth]{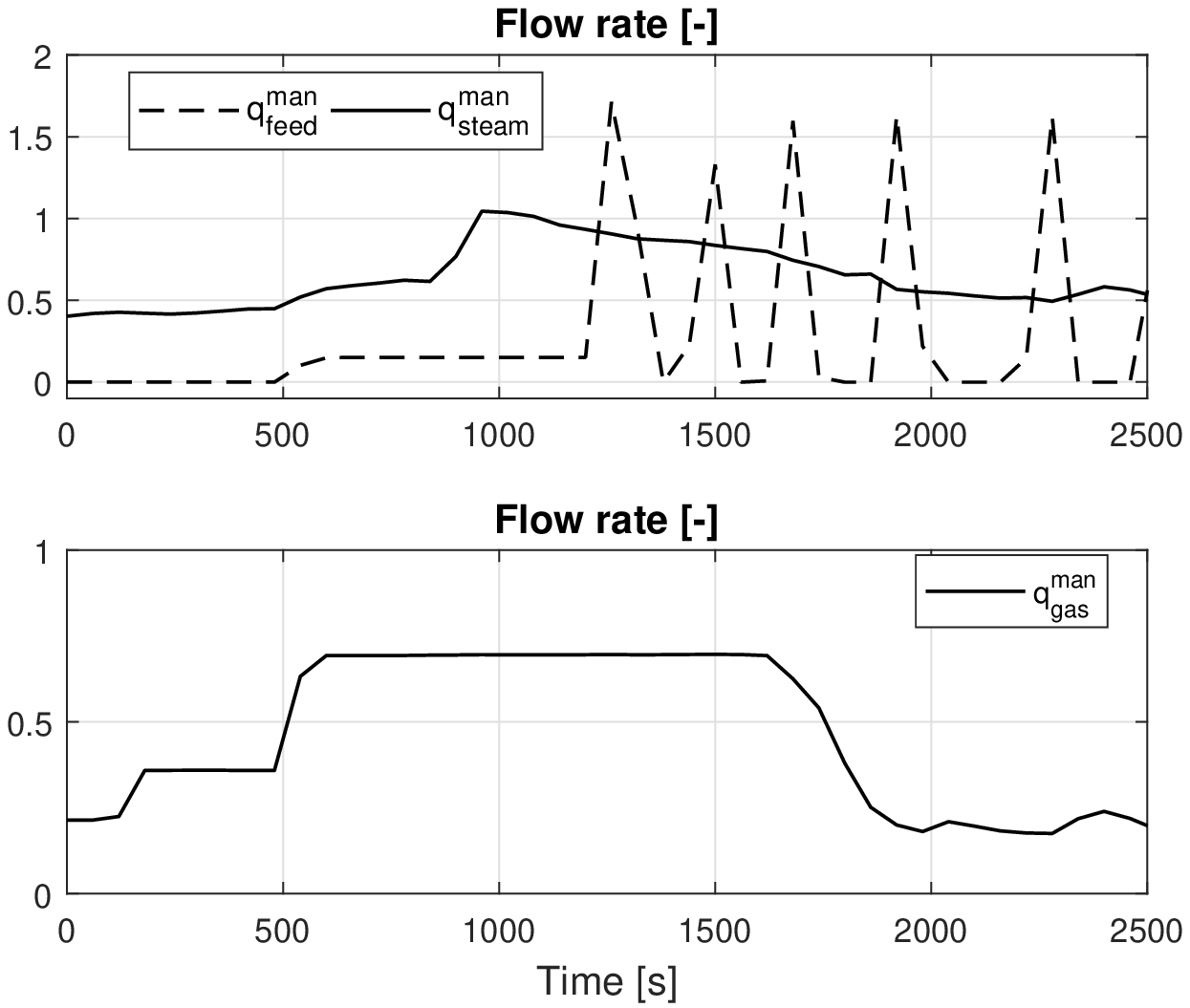}
		\caption{FTB start-up manual procedure - input variables. In the top panel:  steam flow-rate (solid), feedwater flow-rate (dashed). In the bottom panel: gas flow-rate. All graphs are adimensionalized with nominal values, for confidentiality.}
		\label{fig:man_in}
	\end{subfigure}
\caption{Manual FTB start-up }
\end{figure}%
In current industrial practice, start-up procedure is still performed manually. Therefore, a comparison of the proposed LPV-MPC optimal start-up with the typical manual procedure, as extracted by historical data set, is here presented. The human operator during the manual procedure is expected to gradually provide heat to the boiler, at different intermediate levels until a certain water temperature is reached, in order to reduce the thermal stress on the components.
The profiles of the manipulated variables, showing the manual start-up strategy, are displayed in Figure~\ref{fig:man_in}.
The state variables of the manual procedure are shown in Figure~\ref{fig:man_out}: this conservative approach not only induces a longer start-up, as reported in Table~\ref{tab:tab_res}, it even does not guarantee the observance of the thermal stress constraint which is slightly violated, as shown in the second panel of Figure~\ref{fig:man_out}.\\
\begin{figure}[!htb]
	\centering
	\begin{subfigure}[t]{.48\textwidth}
		\centering
		\includegraphics[width=1\linewidth]{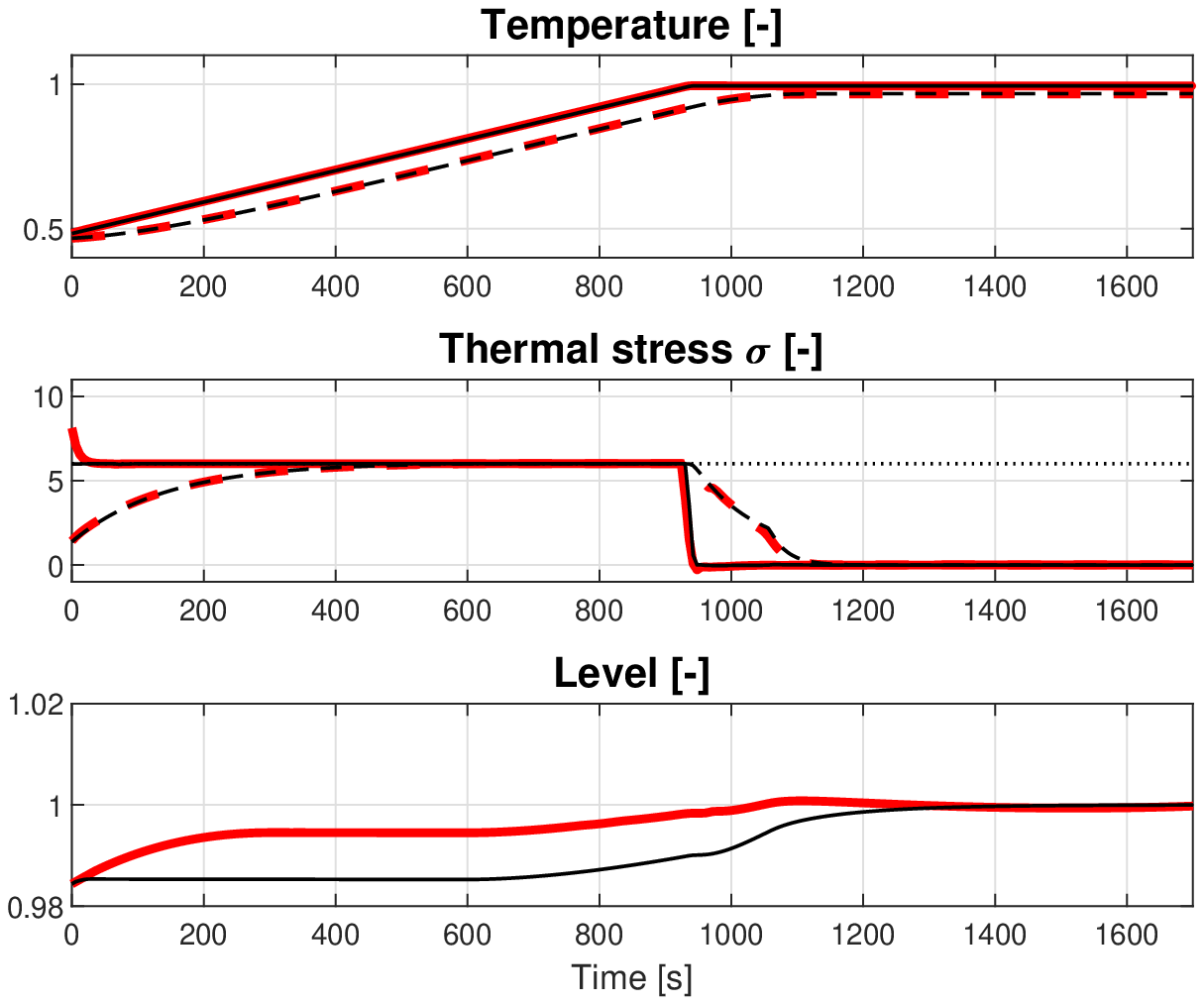}
		\caption{Optimal FTB start-up - output variables. Comparison between open-loop NLP (red bold) and LPV-MPC (black). In the top panel: temperature of fire tubes (solid) and the water temperature (dashed). In the second, the thermal stress on tubes (solid), on the shell (dashed) and its maximum value (dotted). In the third, water level. All graphs are adimensionalized with nominal values, for confidentiality.}
		\label{fig:opt_out}
	\end{subfigure}%
	\hfill
	\begin{subfigure}[t]{0.48\textwidth}
		\centering
		\includegraphics[width=1\linewidth]{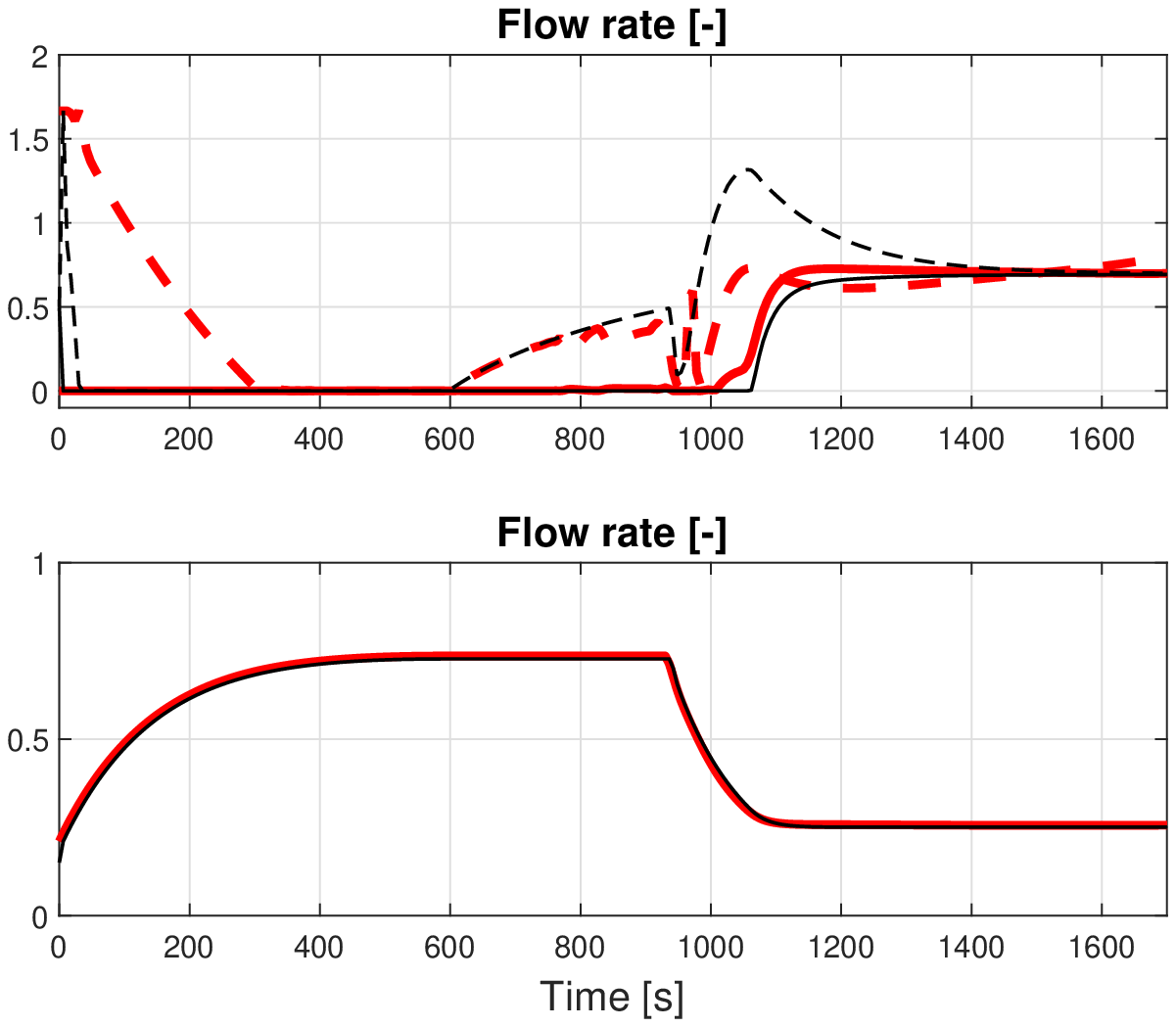}
		\caption{Optimal FTB start-up - input variables. Comparison between open-loop NLP (red bold) and LPV-MPC (black). In the top panel:  steam flow-rate (solid), feedwater flow-rate (dashed). In the bottom panel: gas flow-rate. All graphs are adimensionalized with nominal values, for confidentiality.}
		\label{fig:opt_in}
	\end{subfigure}
\caption{Optimal FTB start-up} 
\end{figure}%
The optimal solution based on LPV-MPC control is, instead, presented in Figure~\ref{fig:opt_out} and Figure~\ref{fig:opt_in}, here compared with the solution of the open-loop NLP: in the former, the output trajectories are shown, while the latter shows the computed optimal inputs.
The open-loop NLP is indeed equivalent to the first step of a closed-loop NMPC, thus providing information about the computational burden required by a NMPC approach, which would apply just the first input of the nonlinear program and then optimize again at the next sampling time, starting from the new initial condition in a receding horizon fashion.\\
\begin{table}[ht]
	\caption{Comparison of start-up procedures: manual start-up (historical data), open-loop NLP (simulation) and closed-loop LPV-MPC (simulation).}
	\label{tab:tab_res}
	\begin{center}
		\begin{tabular}{|l||l||l||l|}
			\hline
			Start-up & Manual & NLP & LPV-MPC\\
			\hline
			type &  & open-loop & closed-loop\\
			\hline
			Duration [s]	($99.3\%$ of SP) 		 & $1678$ & $1026$ & $1032$\\
			\hline
			Energy input [MJ] 			 & \num{1.272e4} & \num{0.947e4} & \num{0.932e4}\\
			\hline
			Thermal stress limit 			 & $\times$ & $\checkmark$ & $\checkmark$\\
			\hline
			Prediction Horizon &			 & $\geq 172$ & $50$\\
			\hline
			Sampling time [s] &			 & $6$ & $6$\\
			\hline
			CPU time [s]  (tot. OCP)&			 &  $\gg 600$ & $0.36$\\
			\hline
		\end{tabular}
	\end{center}
\end{table}%
The overall time required to reach the nominal operating condition is reduced by about 38\% with respect to the conservative manual procedure.
This time reduction is attained by driving quickly the natural gas input closer to the maximum value, while guaranteeing the respect of the constraints. As typical of MPC approaches, the improved performance is obtained by pushing the system closer to the prescribed operating limits.\\
Moreover the LPV-MPC start-up, addressing directly the MIMO system, exploits conveniently the combined effect of the other manipulated variables  in the overall optimization of the start-up trajectory controlling, at the same  time, the water level.\\
As shown qualitatively in Figures~\ref{fig:opt_out}-\ref{fig:opt_in} and quantitatively in Table~\ref{tab:tab_res}, the LPV-MPC method achieves almost identical optimality level with respect to the full horizon nonlinear programming problem. In contrast, the proposed approach offers the main advantages of a reduced CPU time - due to shorter prediction horizon and system linearisation - and closed-loop control of the start-up process.
%
\section{Conclusions} \label{chap:concl}
In this paper we have proposed a hierarchical control structure to optimize the behaviour of a Smart Thermal Energy Grid. The approach is presented on a plant with steam and electric generation units, respectively a fire tube boiler and a combined heat and power internal combustion engine.
The higher hierarchical layer is devoted to optimize the Smart-TEG static behaviour in a 24-hours window, considering  a scenario with price fluctuations and varying demands, by solving an hybrid optimization problem. The solution to the unit commitment problem defines a schedule for the future modes of operation of the systems, as well as the optimal production profiles and the corresponding required combustion gas amounts. Then, at the lower level, a dynamic Model Predictive Controller operates on each individual system to guarantee that the process constraints are fulfilled, with special focus on the boiler. Simulation results are presented to show the effectiveness of the proposed scheme.\\ 
Moreover, the low level scheme is extended with the dynamic control and optimization of the start-up procedure. To address the nonlinearity of the system, we proposed a Linear Parameter-Varying Model Predictive Control approach for the optimization of a nonlinear system subjected to hard constraints. The method exploiting the linearisation of the system along the predicted trajectory is able to address nonlinear system with the advantage of the computational time with respect to NMPC strategy, as approximation of the SQP optimization stopped at the first iteration.\\
The approach is applied to the optimization of the start-up of a nonlinear fire-tube boiler. The simulations showed the performance of the proposed scheme, comparing it to the standard manual approach and a NMPC implementation.\\ 
Future work will address a number of open questions. For example, the dynamic control and optimization of the boiler stand-by phase will be investigated. Even if the existing stand-by controller guarantees that the FTB is maintained in safe operating region, with at the same time a low consumption and a prompt reactivity to return to production modes, the extension of MPC control to this mode might improve the overall efficiency of the system.\\
Moreover, the usage of robust MPC approaches will be analyzed to take into account and manage possible model uncertainties and unmodeled disturbances. \\
Also, we will consider, in the future, the optimization of a large system composed by a number of similar generation units: this will possibly pave the way to tailored distributed optimization and control schemes.%
\section*{Acknowledgements}
This work was partially supported by the SYMBIOPTIMA project. SYMBIOPTIMA has received funding from the European Union's Horizon 2020 research and innovation programme under grant agreement No$^{\circ}$~680426.


\begin{thebibliography}{99}


%
\bibitem{Deloitte2018} Deloitte, (2018). Circular economy in the energy industry, Summary Report
%
\bibitem{SET-Plan2018}  European Technology \& Innovation Platforms Working group 4, (2018). Strategic energy technology plan: increase the resilience and security of the energy system
%
\bibitem{Dukelow2013} Dukelow, S.G., Lipták, B.G., Cheng, X. and Meeker, R.H., (2013). Boiler control and optimization. In Lipták, B.G. (Eds.), Process Control: Instrument Engineers' Handbook, Elsevier Science, pp. 1572-1631
%
\bibitem{Moon2009} Moon, U. C., and Lee, K. Y. (2009). Step-response model development for dynamic matrix control of a drum-type boiler/turbine system. IEEE Transactions on Energy Conversion, 24(2), pp. 423-430.
%
\bibitem{Lu2005} Lu, C.X., Rees, N.W. and Donaldson, S.C., (2005). The use of the {\AA}str{\"o}m-Bell model for the design of drum level controllers in power plant boilers. IFAC Proceedings Volumes, Elsevier, 38(1), pp. 139--144

%
\bibitem{Prasad2002}  Prasad, G., Irwin, G. W., Swidenbank, E. and Hogg, B. W. (2002). A hierarchical physical model-based approach to predictive control of a thermal power plant for efficient plant-wide disturbance rejection. Transactions of the Institute of Measurement and Control, 24(2), 107-128.
%
\bibitem{Pedret2000} Pedret, C., Poncet, A., Stadler, K., Toller, A., Glattfelder, A. H., Bemporad, A. and Morari, M. (2000). Model-varying predictive control of a nonlinear system. Internal Report in Computer Science Dept. ETSE de la Universitat Autònoma de Barcelona.
%
%
\bibitem{FB2005} Feliu-Batlle, V., Rivas Perez, R., Castillo Garcia, F., and Sotomayor Moriano, J. (2005). Fire tube industrial boilers. Fractional control. Automatica and Instrumentation, 365, pp. 90-95.
%
\bibitem{RV2008} Rodriguez-Vasquez, J. R., Perez, R. R., Moriano, J. S., and González, J. P. (2008). Advanced control system of the steam pressure in a fire-tube boiler. IFAC Proceedings Volumes, 41(2), pp. 11028-11033.
%
\bibitem{Hassanein2004} Hassanein, O. I., and Aly, A. A. (2004). Genetic-PID control for a fire tube boiler. In Computational Cybernetics, 2004. ICCC 2004. IEEE pp. 19-24
%
\bibitem{Baldea14} Baldea, M., Harjunkoski, I. (2014). Integrated production scheduling and process control: A systematic review. Computers \& Chemical Engineering, 71, 377-390.
%
\bibitem{Klauco2017}Klau\u{c}o, M., and Kvasnica, M. (2017). Control of a boiler-turbine unit using MPC-based reference governors. Applied Thermal Engineering, 110, pp. 1437-1447.
%
%
\bibitem{Wu2014}Wu, X., Shen, J., Li, Y., and Lee, K. Y. (2014). Hierarchical optimization of boiler–turbine unit using fuzzy stable model predictive control. Control Engineering Practice, 30, 112-123.
%
\bibitem{Baldick1995} Baldick, R. (1995). The generalized unit commitment problem. IEEE Transactions on Power Systems, 10(1), pp. 465-475.

\bibitem{Mitra2013} Mitra, S., Sun, L., and Grossmann, I. E. (2013). Optimal scheduling of industrial combined heat and power plants under time-sensitive electricity prices. Energy, 54, pp. 194-211.

\bibitem{FT2004} Ferrari-Trecate, G., Gallestey, E., Letizia, P., Spedicato, M., Morari, M., and Antoine, M. (2004). Modeling and control of co-generation power plants: a hybrid system approach. IEEE Transactions on Control Systems Technology, 12(5), pp. 694-705.
%
\bibitem{Ashok2002} Ashok, S., Banerjee, R. (2003). Optimal operation of industrial cogeneration for load management. IEEE Transactions on power systems, 18(2), 931-937.
%
\bibitem{Hawkes2009} Hawkes, A. D.,  Leach, M. A. (2009). Modelling high level system design and unit commitment for a microgrid. Applied energy, 86(7-8), 1253-1265.
%
\bibitem {Anand2019} Anand, H., Narang, N., Dhillon, J. S. (2019). Multi-objective combined heat and power unit commitment using particle swarm optimization. Energy, 172, 794-807.
%
%
%
%
%






\bibitem{Hubel2017} Hübel, M., Meinke, S., Andrén, M. T., Wedding, C., Nocke, J., Gierow, C., Hassel E., Funkquist, J. (2017). Modelling and simulation of a coal-fired power plant for start-up optimisation. Applied Energy, 208, 319-331.

\bibitem{Sindareh2014} Sindareh-Esfahani, P., Habibi-Siyahposh, E., Saffar-Avval, M., Ghaffari, A., Bakhtiari-Nejad, F. (2014). Cold start-up condition model for heat recovery steam generators. Applied Thermal Engineering, 65(1-2), 502-512.

\bibitem{Albanesi2006} Albanesi, C., Bossi, M., Magni, L., Paderno, J., Pretolani, F., Kuehl, P., Diehl, M. (2006). Optimization of the start-up procedure of a combined cycle power plant. In Decision and Control, 2006 45th IEEE Conference, pp. 1840-1845.

\bibitem{Casella2011} Casella, F., Farina, M., Righetti, F., Scattolini, R., Faille, D., Davelaar, F.,Tica A. , Gueguen H.  Dumur, D. (2011). An optimization procedure of the start-up of combined cycle power plants. In 18th IFAC World Congress pp. 7043-7048.


\bibitem{Kru2004} Krüger, K., Franke, R., Rode, M. (2004). Optimization of boiler start-up using a nonlinear boiler model and hard constraints. Energy, 29(12-15), 2239-2251.

\bibitem{Franke2003} Franke, R., Rode, M., Krüger, K. (2003). On-line optimization of drum boiler startup. In Proceedings of the 3rd International Modelica Conference

\bibitem{Belkhir2015} Belkhir, F., Cabo, D. K., Feigner, F., Frey, G. (2015). Optimal startup control of a steam power plant using the JModelica platform. IFAC-PapersOnLine, 48(1), 204-209.

\bibitem{Taler2015} Taler, J., Wlglowski, B., Taler, D., Sobota, T., Dzierwa, P., Trojan, M., Madejski P., Pilarczyk, M. (2015). Determination of start-up curves for a boiler with natural circulation based on the analysis of stress distribution in critical pressure components. Energy, 92, 153-159.

\bibitem{Taler2015_2} Taler, J., Dzierwa, P., Taler, D., Harchut, P. (2015). Optimization of the boiler start-up taking into account thermal stresses. Energy, 92, 160-170.

\bibitem{Frank2009} Franke, R. Vogelbacher, L. (2009). Nonlinear model predictive control for cost optimal startup of steam power plants (Nichtlineare modellprädiktive Regelung zum kostenoptimalen Anfahren von Dampfkraftwerken). Automatisierungstechnik, 54(12), pp. 630-637.
%
\bibitem{Cominesi2017}  {Raimondi Cominesi}, S., Farina, M., Giulioni, L., Picasso, B., and Scattolini, R. (2017). A two-layer stochastic model predictive control scheme for microgrids. IEEE Transactions on Control Systems Technology, 1-13.
%
\bibitem{Limon2008} Limón, D., Alvarado, I., Alamo, T., Camacho, E. F. (2008). MPC for tracking piecewise constant references for constrained linear systems. Automatica, 44(9), pp. 2382-2387.




\bibitem{Spinelli2018} Spinelli, S., Farina, M., Ballarino, A. (2018). A hierarchical optimization-based scheme for combined Fire-tube Boiler/CHP generation units. Proceedings of European Control Conference 2018, pp. 416-421.

\bibitem{Spinelli2019} Spinelli, S., Farina, M., Ballarino, A. (2019). An optimal control of start-up for nonlinear fire-tube boilers with thermal stress constraints, Proceedings of European Control Conference 2019, pp. 2362-2367.

%
\bibitem{Lygeros2008}  Lygeros, J., Tomlin, C., and Sastry, S. (2008). Hybrid systems: modeling, analysis and control.
\bibitem{Guzzella2004}  L. Guzzella and C. H. Onder,(2004).  Introduction to modeling and control of internal combustion engine systems. Springer
%
\bibitem{Videla2007}Videla, J. I., and Lie, B. (2007). State/parameter estimation of a small-scale CHP model. 48th Scandinavian Conference on Simulation and Modeling (SIMS 2007); 30-31; Goteborg. Linköping University Electronic Press.



\bibitem{Ortiz2011} Ortiz, F. G. (2011). Modeling of fire-tube boilers. Applied Thermal Engineering, 31(16), pp. 3463-3478.

\bibitem{Astrom2000} Astrom, K.J., Bell, R.D. (2000). Drum-boiler dynamics. Automatica 36  pp. 363-78.

\bibitem{Bejan2003} Bejan, A., and Kraus, A. D. (2003). Heat transfer handbook (Vol. 1). John Wiley \& Sons.

\bibitem{Bemporad1999} Bemporad, A., and Morari, M. (1999). Control of systems integrating logic, dynamics, and constraints. Automatica, 35(3), pp. 407-427. 
%

\bibitem{SteamTab} Wagner, W. and Kretzschmar, H. J. (2007). International steam tables-properties of water and steam based on the industrial formulation IAPWS-IF97. Springer Science \& Business Media.

\bibitem{EN12952} DIN-EN (2008). 12952-3: Water-tube boilers and auxiliary installations-Part 3: Design and calculation for pressure parts.

\bibitem{SteamTool}  Mikofski M. (2019). IAPWS-IF97 (https://www.github.com/mikofski/IAPWS\textunderscore IF97), GitHub. Retrieved May 18, 2019. 

\bibitem{HYS2004} Torrisi, F. D., and Bemporad, A. (2004). HYSDEL -  a tool for generating computational hybrid models for analysis and synthesis problems. IEEE transactions on control systems technology, 12(2), pp. 235-249.

\bibitem{Simon2006} Simon, D. (2006). Optimal state estimation: Kalman, H infinity, and nonlinear approaches. John Wiley \& Sons.
%
\bibitem{Mitra94} Mitra, G., Lucas, C., Moody, S., and Hadjiconstantinou, E. (1994). Tools for reformulating logical forms into zero-one mixed integer programs. European Journal of Operational Research, 72(2), pp. 262-276.
%
\bibitem{Williams2013} Williams, H. P. (2013). Model building in mathematical programming. John Wiley \& Sons, pp. 154-165.

\bibitem{Falcone08} Falcone, P., Borrelli, F., Tseng, H. E., Asgari J. and Hrovat D. (2008). Linear time-varying model predictive control and its application to active steering systems: stability analysis and experimental validation. International Journal of Robust and Nonlinear Control, 18, 862-875.


\bibitem{Diehl2005} Diehl, M., Bock, H. G., Schlöder, J. P. (2005). A real-time iteration scheme for nonlinear optimization in optimal feedback control. SIAM Journal on control and optimization, 43(5), 1714-1736.




\end{thebibliography}
\end{document}